\documentstyle[amscd,12pt,psamsfonts,leqno,amssymb]{amsart}
\input xypic

\newtheorem{defn}{Definition}
\newtheorem{thm}[defn]{Theorem}

\newtheorem{cor}[defn]{Corollary}
\newtheorem{lem}[defn]{Lemma}
\newtheorem{prop}[defn]{Proposition}

\theoremstyle{remark}
\newtheorem{rem}{Remark}
\theoremstyle{remark}

\numberwithin{equation}{section} \numberwithin{defn}{section}

\renewcommand\sp{\operatorname{Spec}}
\renewcommand\sf{\operatorname{Spf}}
\newcommand\grv{{\operatorname{Gr}}(V)}
\newcommand\gr{\operatorname{Gr}}
\newcommand\pgrv{{\operatorname{Gr}_{\text{{\rm I}}}}(V)}
\newcommand\nm{\operatorname{Nm}}

\renewcommand\hom{\operatorname{Hom}}
\renewcommand\det{\operatorname{det}}
\newcommand\Det{\operatorname{Det}}

\renewcommand\ker{\operatorname{Ker}}
\newcommand\res{\operatorname{Res}}
\newcommand\limi{\varinjlim}
\newcommand\limil[1]{\underset{#1}\varinjlim\,}
\newcommand\limp{\varprojlim}
\newcommand\limpl[1]{\underset{#1}\varprojlim\,}
\newcommand\proj{\operatorname{Proj}}
\newcommand\aut{\operatorname{Aut}}
\renewcommand\hom{\operatorname{Hom}}

\newcommand\tr{\operatorname{Tr}}

\renewcommand\deg{\operatorname{deg}}
\newcommand\id{\operatorname{Id}}

\newcommand\Rad{\operatorname{Rad}}

\newcommand\pic{\operatorname{Pic}}

\newcommand\prym{\operatorname{Prym}}
\newcommand\fprym{\underline\prym}
\newcommand\fm{\underline\M}

\renewcommand\o{{\mathcal O}}
\renewcommand\j{{\mathcal J}}

\newcommand\p{{\mathcal P}}
\newcommand\G{{\mathbb G}}
\renewcommand\L{{\mathcal L}}

\newcommand\M{{\mathcal M}}

\newcommand\U{{\mathcal U}}

\newcommand\Z{{\mathbb Z}}

\newcommand\C{{\mathbb C}}

\newcommand\w{\widehat}
\newcommand\wtilde{\widetilde}

\newcommand\beq{
                     \setcounter{equation}{\value{defn}}\addtocounter{defn}1
                     \begin{equation}}

\begin{document}

\title[Prym varieties and curves with automorphisms]
{Prym varieties, \\
 curves with automorphisms \\
 and the Sato
Grassmannian}

\author[E. G\'omez, J. M. Mu\~noz  and F. J. Plaza]{E. G\'omez
Gonz\'alez$\;^{(1)}$ \\ J. M. Mu\~noz Porras$\;^{(1)}$ \\  F. J. Plaza
Mart\'{\i}n$\;^{(2)}$}

\address{$\;^{(1)}$ Departamento de Matem\'aticas, Universidad de
Salamanca,  Plaza
        de la Merced 1-4
        \\
        37008 Salamanca. Spain.
        \\
         Tel: +34 923294460. Fax: +34 923294583}

\address{$\;^{(2)}$ Departamento de Matem\'aticas, Universidad
Aut\'onoma de Madrid,
         Cantoblanco
        \\
        28049 Madrid. Spain.
        \\
         Tel: +34 913967635. Fax: +34 913974889}

\thanks{
\indent This work is partially supported by the research contracts
BFM2000-1327 and BFM2000-1315 of DGI and  SA064/01 of JCyL. The third
author is also supported by MCYT ``Ram\'on y Cajal'' program. \\
\indent {\it E-mail addresses}: esteban@@usal.es, jmp@@usal.es,
francisco.plaza@@uam.es}
\email{esteban@@usal.es} \email{jmp@@usal.es}
\email{francisco.plaza@@uam.es}

\begin{abstract}
The aim of the paper is twofold. First, some results of Shiota
and Plaza-Mart\'{\i}n on Prym varieties of curves with an involution
are generalized to the general case of  an arbitrary automorphism
of prime order. Second,  the equations defining the moduli space
of curves with an automorphism of prime order  as a subscheme of
the Sato Grassmannian are given.
\end{abstract}


\maketitle



\section{Introduction}

The main objectives of this paper are as follows. First, to
extend some results of T. Shiota (\cite{Sh2}) and F. J.
Plaza-Mart\'{\i}n (\cite{PlPrym}) on Prym varieties of curves with an
involution to the general case of Prym varieties associated with
curves with an arbitrary automorphism of prime order. Second, to
give an explicit description of the equations defining the moduli
space of curves with an automorphism of prime order (together with
a formal trivialization at some points) as a subscheme of the Sato
Grassmannian.

The results of this paper are intended to be a first step to solving
a Schottky-type problem for curves with automorphisms. The basic
concepts used in our approach are the formal spectral covers and
their formal Jacobians, which were introduced in \cite{MP3} when
studying the Hurwitz schemes.

We define the notion of formal Prym variety associated with a formal
Galois cover and prove that this formal Prym satisfies all the
expected properties (i.e. the analogues of those of classical Prym
varieties, see Proposition~\ref{prop:expectedpropertiesprym}).
This formal Prym variety acts on the infinite Grassmannian of $V$,
where $V$ is the $\C((z))$-algebra of the formal cover.

In~\S3, we study two subschemes of the infinite Grassmannian of
$V$ that are relevant for our purposes. The first is the set
of  points of the Grassmannian invariant under an automorphism of
$V$; that is, subspaces $U\subset V$ stable under the
automorphism. The second one consists of a generalization of the
Grassmannian of maximal totally isotropic subspaces considered in
~\cite{PlPrym,Sh2} when  studying the BKP hierarchy. In the
present setting, the bilinear pairing introduced in those papers
is replaced by a multilinear form on $V$.

The Krichever map for Pryms is studied in~\S4 in order to obtain a
characterization of Prym varieties
(Theorem~\ref{thm:characterizationofpryms}) as subvarieties of
the Sato Grasmmannian. This statement (together with
Theorem~\ref{thm:algebraizableiffgeometric}) is to be thought of as
an analogue of Mulase's characterization of Jacobians (\cite{Mul})
and of Li-Mulase of Pryms (\cite{LiMulasePrym}) since it is
expressed in terms of the finiteness of the orbits (modulo
$\bar\Pi^+$) of the action of the formal group $\Pi$ on the
Grassmannian.

Finally, in~\S5 we state the explicit equations of the moduli
spaces of curves with automorphisms as subvarieties of the
infinite Grassmannians. To this end, we use some ideas and results
on the equations of the moduli spaces of curves that were given by
two of the authors in a previous paper (\cite{MP}).

Through this paper we will assume that the base field is $\C$, the
field of complex numbers.

The authors wish to express their gratitude to referee for his/her
valuable comments.


\section{Formal Geometry of Coverings}

This section is concerned with the generalization of the notions
of the formal curve, the formal jacobian, etc. (see~\cite{AMP,MP3}) for
the case of a covering of curves of group $\Z/p$.

     Henceforth,  $V$ will be a $\C((z))$-algebra with an action of
the group $\Z/p$, where $p$ is a prime number such that
$V^{\Z/p}=\C((z))$. We shall choose a generator $\sigma$ of
$\aut_{\C((z))}V$. We consider the following two cases:
\begin{enumerate}
\item[(a)] {\sl Ramified case:}
$V=\C((z_1))$ where the $\C((z))$-algebra structure is given by
mapping $z$ to $z_1^p$ and $\sigma(z_1)=\xi z_1$ ($\xi$ being a
primitive $p$-th root of $1$ in $\C$). We will set $V^+=\C[[z_1]]$ and
$V^-=z_1^{-1}\C[z_1^{-1}]$.
\item[(b)] {\sl Non-ramified case:} $V=\C((z_1))\times \dots\times
\C((z_p))$ where the
$\C((z))$-algebra structure is given by mapping $z$ to
$(z_1,\dots, z_p)$ and $\sigma(z_i)=z_{i+1}$ (for $i<p$) and
$\sigma(z_p)=z_1$. We set $V^+=\C[[z_1]]\times
\dots\times \C[[z_p]]$ and $V^-=z_1^{-1}\C[z_1^{-1}]\times \dots\times
z_p^{-1}\C[z_p^{-1}]$.
\end{enumerate}
     In both cases, we have {\it distinguished} bases
of $V$ as a $\C((z))$-vector space; namely,
$\{1,z_1,\dots,z_1^{p-1}\}$ in the first case, and
$\{z_1,\dots,z_p\}$ in the second one. Observe that
$\sigma(V^+)=V^+$ and that $\sigma(V^-)=V^-$.

\begin{rem}
If $p$ is not a prime number a more general $\C((z))$-algebra $V$
must be considered. Then, the results of this paper could be
generalized to the case of a cyclic group of automorphisms.
However, for the sake of clarity, we will restrict ourselves to
the case of $p$ being a prime number.
\end{rem}

Following~\cite{MP3}, Definition~2.1, we define the formal base
curve as the formal $\C$-scheme $\w C:=\sf \C[[z]]$ and the formal
spectral cover as the formal $\C$-scheme $\w C_V:=\sf V^+$.

Let $\Gamma$ denote the formal group scheme representing the
functor: {\small $$
\begin{aligned}
\left\{{\scriptsize\begin{gathered}\text{category of}\\
\text{formal $\C$-schemes}\end{gathered}}\right\}\,&
\rightsquigarrow\,\left\{{\scriptsize\begin{gathered}\text{category}\\\text{of
groups}\end{gathered}}\right\} \\
S\,& \rightsquigarrow\,  \big(\C((z))\hat{\underset{\C}\otimes}
H^0(S,\o_S)\big)^*_0 \,=\, \left\{{\scriptsize
\begin{gathered} \text{invertible elements}\\
\text{of }\C((z))\hat{\underset{\C}\otimes}  H^0(S,\o_S)
\end{gathered}}\right\}_0
\end{aligned}
$$}
where the subscript $0$ denotes the connected component of the
identity. Let $\Gamma^+$  be the group
representing the subfunctor of $\Gamma$, consisting of those
elements lying in $\C[[z]]\hat\otimes H^0(S,\o_S)$, and let the
      jacobian of the formal curve,  $\j(\w C)$,  be the group
representing the subfunctor of $\Gamma$, consisting of those
elements lying in
$\C[z^{-1}]\hat\otimes H^0(S,\o_S)$ whose constant term is
equal to $1$. Observe that its
$R$-valued points ($R$ being a $\C$-algebra) is the set:
      {\small $$\left\{
      \begin{gathered} \text{polynomials $1+\sum_{j=-n}^{-1}a_j z^j$
      of $R[z^{-1}]$} \\
      \text{where $a_j\in\Rad(R)$ for all $j$}
      \end{gathered}\right\}$$}
Note that the multiplicative group $\G_m$ is
contained in $\Gamma$ and that
$\Gamma=\j(\w C)\times\G_m\times\Gamma^+$.

In the ramified case, to define $\Gamma_V$ (resp. $\Gamma^+_V$)
one replaces $\C((z))$ (resp. $\C[[z]]$) by $\C((z_1))$ (resp.
$\C[[z_1]]$). The formal jacobian of the formal spectral cover
$\j(\w C_V)$ is defined by replacing $\C[z^{-1}]$ by
$\C[z_1^{-1}]$ and it turns out that $\Gamma_V=\j(\w
C_V)\times\G_m\times \Gamma_V^+$.

In the non-ramified case, $\Gamma_V$ (resp. $\Gamma^+_V$) is
defined by replacing $\C((z))$ (resp. $\C[[z]]$) by
$\C((z_1))\times \dots\times \C((z_p))$ (resp.
$\C[[z_1]]\times\dots\times \C[[z_p]]$). To define the formal
jacobian of the formal spectral cover, $\j(\w C_V)$, one replaces
$\C[z^{-1}]$ by $\C[z_1^{-1}]\times \dots\times\C[z_p^{-1}]$; one then
has that $\Gamma_V=\j(\w
C_V)\times(\G_m\times\overset{p}\dots\times\G_m)\times
\Gamma_V^+$.

\begin{rem}\label{rem:jisdirectlim}
It is worth noticing that the formal scheme $\Gamma$ is a direct
limit of affine group schemes $\Gamma^n$. For $n\geq 0$, one
considers: { $$\Gamma^n\,:=\, \left\{{\small
\begin{gathered} \text{invertible elements}\\
\text{of }z^{-n}\C[[z]]\hat{\underset{\C}\otimes}  H^0(S,\o_S)
\end{gathered}}\right\}_0$$}
and observes that it is the spectrum of $\C[u_n,\dots,
u_1,v_0,v_1,\dots]/I_n$, where $u_i$ has degree $i$ and $I_n$ is
the ideal generated by the monomials on $u$ of total degree
greater than $n$.  The results of~\cite{EGA},~I.\S10.6 imply that
$\Gamma=\limi \Gamma^n$.

Analogous arguments show that $\j(\w C)$, $\Gamma_V$ and $\j(\w
C_V)$ are limits of affine group schemes.
\end{rem}

Standard calculation shows that $\j(\w C_V)$ is the formal
spectrum of the ring: {\small $$\C\{\{t_1,t_2,\dots\}\}$$}
\noindent ($\C\{\{t_1,t_2,\dots\}\}$ denotes the inverse limit
$\limpl{n}\C[[t_1,\dots,t_n]]$) for the ramified case and: {\small
$$\C\{\{t_1^{(1)},t_2^{(1)},\dots\}\}\w\otimes \dots \w\otimes
\C\{\{t_1^{(p)},t_2^{(p)},\dots\}\}$$} for the non-ramified one.

Since $\sigma$ acts on $V^+$ as a homomorphism of
$\C[[z]]$-algebras, it gives rise to an order $p$ automorphism:
$$\sigma\colon \w C_V\,\overset{\sim}\longrightarrow\, \w C_V$$
such that the quotient $\w C_V/<\!\sigma\!>=\sf
(V^+)^{<\!\sigma\!>}=\sf \C[[z]]$ is the formal base curve.

     From the explicit expression of the elements of $\j(\w C_V)$, it is
easy to check that $\sigma$ induces an automorphism
$\sigma^*:\j(\w C_V)\overset{\sim}\to \j(\w C_V)$. The expression
in terms of rings is:
$$ t_j\,\longmapsto \, \sigma^*(t_j)=\xi^{-j} t_j$$
and
$$t_j^{(i)}\,\longmapsto \, \sigma^*(t_j^{(i)})\,=\,
\cases t_j^{(i-1)} & \quad i>1 \\ t_j^{(p)} & \quad i=1
\endcases$$ respectively.

The Abel morphism of degree $1$ is the morphism:
$$\w C_V\,\longrightarrow\, \j(\w C_V)$$
defined by the series:
{\small
$$\exp(\sum_{j>0}\frac{\bar z_1^j}{j z_1^j})$$}
in the ramified case (here $\o_{\w C_V}$ is identified with
$\C[[\bar z_1]]$) and by the $p$-uple of series: {\small
$$\big(\exp(\sum_{j>0}\frac{\bar z_1^j}{j z_1^j}), \dots,
\exp(\sum_{j>0}\frac{\bar z_p^j}{j z_p^j})\big)$$} in the
non-ramified one ($\o_{\w C_V}$ is identified with $ \C[[\bar
z_1]]\times\dots\times \C[[\bar z_p]]$).

$V$ being a finite $\C((z))$-algebra, we have the norm, which is a
group homomorphism:
$$\nm \colon V^*\,\longrightarrow\, \C((z))^*$$
that maps an element $v\in V^*$ to the determinant of the
homothety defined by itself or, equivalently, the product of its
transforms under $\sigma^i$ ($0\leq i\leq p-1$).
Recalling the definition of formal jacobians of $\w C_V$
and $\w C$, one has that the
norm gives rise to a morphism of formal groups schemes:
$$\nm\colon \j(\w C_V) \,\longrightarrow\, \j(\w C)$$
which corresponds to the ring homomorphism:
{\small $$
\begin{aligned}
\C\{\{\bar t_1,\bar t_2,\dots\}\}& \,\longrightarrow\,
\C\{\{t_1,t_2,\dots\}\}  \\
\bar t_i &\,\longmapsto \, \nm(\bar t_i)= p\cdot  t_{ip}
\end{aligned}$$}
in the ramified case and:
{\small $$
\begin{aligned}
\C\{\{\bar t_1,\bar t_2,\dots\}\}& \,\longrightarrow\,
\C\{\{t_1^{(1)},\dots\}\}\w\otimes\dots \w\otimes
\C\{\{t_1^{(p)},\dots\}\} \\
\bar t_i &\,\longmapsto \,\nm(\bar t_i)= t_i^{(1)}+\dots +
t_i^{(p)}
\end{aligned}$$}
in the non-ramified case.

Finally, since the quotient map $\pi:\w C_V\to \w C$ corresponds
to the canonical inclusion $\C[[z]]\hookrightarrow V^+$, it also
induces a group morphism:
$$\pi^* \colon \j(\w C)\,\hookrightarrow\, \j(\w C_V)$$
whose expression in terms of rings is:
$$t_j\,\longmapsto
\, \pi^*(t_j)\,=\, \cases  \bar t_{j/p} &\quad\text{if }j=\dot{p} \\ 0
&\quad\text{otherwise}\endcases$$
for the ramified case and:
$$t_j^{(i)}\,\longmapsto \, \pi^*(t_j^{(i)})=\bar t_j$$
for the
non-ramified one.

Note that the set of invariant elements of the formal jacobian of
the spectral cover, $\j(\w C_V)^{<\!\sigma^*\!>}$, is the formal
jacobian of the formal base curve, $\j(\w C)$.

\subsection{Formal Prym varieties}\hfill

Let us introduce two  important subgroups of $\Gamma_V$. Note that
the norm can be generalized to the relative case since
$V\w\otimes\o_S$ is a finite $\C((z))\w\otimes\o_S$-algebra.

\begin{defn}
The formal group $\Pi$ is the formal group scheme representing the
subfunctor of $\Gamma_V$ given by:
$$S\,\rightsquigarrow\,\Pi(S)\,:=\, \{ g\in\Gamma_V(S)\text{ such that
}\nm(g)\in\o_S\}$$ and $\Pi^+$ is defined by $\Pi\cap \Gamma_V^+$.

The formal Prym variety of $(V,V^+,\sigma)$ is the formal group
scheme   given by:
$$\p(\w C_V)\,:=\ker(\nm\vert_{\j(\w C_V)})\,=\,
\{g\in \j(\w C_V)\;\vert\; \nm(g)=1\}$$
\end{defn}

Let us state some of the properties of these subgroups of $\Gamma_V$.

\begin{prop}\label{prop:expectedpropertiesprym}
The following properties hold:
\begin{enumerate}
\item
$\p(\w C_V)= \Pi\cap  \j(\w C_V)$;
\item
the map $\sigma^*-\id:\Gamma_V\to \Gamma_V$ induces a surjection
$\j(\w C_V)\to\p(\w C_V)$. Further, it gives a map
$\sigma^*-\id:\w C_V\to \p(\w C_V)$;
\item
the multiplication of series gives rise to an isomorphism: {\small
$$m:\p(\w C_V)\times \j(\w C)\to \j(\w C_V)\ ;$$}
\item
there is an isomorphism: {\small $$(\sigma^*-\id,\nm):\j(\w
C_V)\to \p(\w C_V)\times \j(\w C)\ ;$$}
\item
the compositions $(\nm,\sigma^*-\id)\circ m$ and $m\circ
(\nm,\sigma^*-\id)$ map an element to its $p$-th power.
\end{enumerate}
\end{prop}

\begin{pf}
Note that the set of $R$-valued points of $\p(\w C_V)$ is the
following subset of $\j(\w C_V)(R)$:
{\small
$$\big\{\exp(\sum_{j} t_j z_1^{-j})\,\in \, \j(\w
C_V)(R)\text{ such that $t_j=0$ for }j= \dot{p}  \Big\}$$}
for the ramified
case and:
{\small
$$\Big\{\big(\exp(\sum_{j>0} t^{(1)}_j z_1^{-j}),\dots
,\exp(\sum_{j>0} t^{(p)}_j z_p^{-j})\big) \text{ such that $
\sum_{i=1}^p t_j^{(i)}=0$ for }j>0\Big\}$$} for the
non-ramified one.

Now, the claims follow easily from the explicit expressions of the
morphisms $\nm$, $\pi^*$ and $\sigma^*$ given in the previous
section.
\end{pf}

\section{Vector-valued infinite Grassmannians}

We begin this section by reviewing from~\cite{AMP,MP3} some
definitions and main properties of infinite Grassmannians.

From~\cite{AMP} we learn
that the Grassmannian functor of $(V,V^+)$ is representable by a
$\C$-scheme whose set of rational points is:
  {\small $$\left\{\begin{gathered}
  \text{subspaces $U\subset V$ such that }U\to V/V^+
  \\
  \text{ has finite dimensional kernel and cokernel}
  \end{gathered}
  \right\}$$}
The connected components of this scheme are indexed by the
Poincar\'e-Euler characteristic of $U\to V/V^+$. The connected
component of index $m$ will be denoted by $\gr^m(V)$. This scheme
is  equipped with the determinant bundle, which is the determinant
of the complex of $\o_{\grv}$-modules:
$$\L\,\longrightarrow\, V/V^+$$
where $\L$ is the universal submodule of $\grv$ and $\L\to V/V^+$
is the natural projection.

Recall that $\grv$ is endowed with an action of the group
$\Gamma_V$:
$$\Gamma_V\times\grv\,\longrightarrow\,\grv$$
since the elements of $\Gamma_V$ act by homotheties on $V$.

The Baker-Akhiezer function  of a point $U\in\grv$ has been
introduced in~\S3 of~\cite{MP3}. This is a function on $\w
C_V\times\j(\w C_V)$ that will be denoted by
$\psi_U(z_\centerdot,t)$, where $z_\centerdot$ is $z_1$ (resp.
$(z_1,\dots, z_p)$) in the ramified case (resp. in the
non-ramified case). It is appropriate to point out that when
defining the Baker-Akhiezer function of $U\in \gr^m(V)$, one has
to choose an element $v_m\in V$ such that the kernel and cokernel
of $U\to V/v_mV^+$ will have the same dimension.

Following that paper, we choose $v_m$ to be  $z_1^m$ in the
ramified case and  $(z_1^{q+1},\dots, z_r^{q+1}, z_{r+1}^q,\dots,
z_p^q)$ (where $m=pq+r$) in the non-ramified one.

The main property of Baker-Akhiezer functions is  that
$\frac{v_m}{z_\centerdot}\psi_U(z_\centerdot,t)$ can be
understood as a generating function of the vectors of
$U\in\gr^m(V)$ as a subspace of $V$. That is, every $u\in U$ is
obtained by evaluating
$\frac{v_m}{z_\centerdot}\psi_U(z_\centerdot,t)$ at certain values
of $t$, and, conversely, that this expression belongs to $U$ for all
values of $t$.

Since $V$ is a finite $\C((z))$-algebra, there is a natural linear
map:
$$\tr\colon V\,\longrightarrow\,\C((z))$$
which assigns to an element $v\in V$ the trace of the homothety of
$V$ defined by $v$. Furthermore, it gives rise to a
non-degenerated metric on $V$ as $\C((z))$-vector space, which
will be denoted by $\tr$ again.

This metric allows us to consider the following pairing:
\begin{equation}\label{eq:restrace}
\begin{aligned}
   V\times V &
\,\longrightarrow\, \C \\
(a,b) & \,\longmapsto \res_{z=0}\tr(a,b) dz
\end{aligned}
\end{equation}

Since this pairing is non-degenerate, we have an involution of
$\grv$ which  maps a point $U$ to its orthogonal $U^\perp$. This
involution  sends the connected component $\gr^m(V)$ to
$\gr^{1-m-p}(V)$ in the ramified case and to $\gr^{-m}(V)$ in the
non-ramified one.

Finally, the adjoint Baker-Akhiezer function of $U$ is defined by:
$$\psi_U^*(z,t)\,:=\,\psi_{U^\perp}(z,-t)$$

\subsection{Stable points of $\grv$ under an automorphism}\hfill

We shall now take into account the additional structure provided
by the automorphism $\sigma$. We shall show that the subset of
points of $\grv$ invariant under the automorphism is a closed
subscheme and shall compute how the Baker-Akhiezer function
behaves. In fact, we shall consider a more general setting.

Let us denote by $\aut^+_{\C((z))}V$ the group of
$\C((z))$-algebra automorphisms of $V$ that leave  $V^+$ stable. Let
$G$ be a finite group and $\rho:G\hookrightarrow
\aut^+_{\C((z))}V$ be a linear representation of $G$.

Observe that if $A\subset V$ is a subspace commensurable with
$V^+$, then  $\rho(\sigma)A$ is also commensurable with $V^+$ for
all $\sigma\in G$. Therefore, the representation $\rho$ yields an
action on $\grv$:
\begin{equation}\label{eq:actiongroupongrv}
\begin{aligned}
G\times\grv&\,\longrightarrow\, \grv\\
(\sigma,U)&\,  \longmapsto \rho(\sigma)U
\end{aligned}
\end{equation}

Let us set:
$$\grv^G\,:=\,\{U\in\grv\text{ such that }\rho(\sigma)U=U\text{ for all
}\sigma\in G\}$$ Then, one has that:
$$\grv^G\,=\,\underset{\sigma\in G}\bigcap \grv^\sigma$$
where:
$$\grv^\sigma\,:=\,\{U\in\grv\text{ such that }\rho(\sigma)U=U\}$$
Note that  $\grv^\sigma$ is a closed subscheme of $\grv$
because $\grv$ is separated (\cite{AMP}, Theorem~2.15) and,
therefore, $\grv^G$ is also a closed subscheme.

Since the Baker-Akhiezer and adjoint Baker-Akhiezer functions
might be thought of (up to a factor) as generating functions, we
obtain the following characterization of $\grv^\sigma$ for
$\sigma^p=\id$:

\begin{thm}\label{thm:characterizationgrsigma}
Let $\sigma\in \aut^+_{\C((z))}V$ with $\sigma^p=\id$. A closed
point $U\in\grv$ is a point of $\grv^\sigma$ if and only if its
Baker-Akhiezer function satisfies the equation:
$$
\res_{z=0}\tr\Big(\frac1{z_1}\psi_{\rho(\sigma)U}(z_1,t)
\psi_U^*(z_1,s)\Big)\frac{dz}z\,=\,0$$ in the ramified case, and:
$$
\res_{z=0}\tr\Big(\frac1{{z_\centerdot^2}}\psi_{\rho(\sigma)U}({z_\centerdot},t)
\psi_U^*({z_\centerdot},s)\Big)dz\,=\,0$$ in the non-ramified one.
\end{thm}

\begin{pf}
   Let us begin with the ramified case.
Given a point $U\in\gr^m(V)$, we know that $U^\perp$ belongs to
$\gr^{1-m-p}(V)$. In particular,
   we know that $\frac{v_m}{z_1}\psi_U(z_1,t)$ (resp.
$\frac{v_{1-m-p}}{z_1}\psi^*_U(z_1,t)$) is a generating function
of $U$ (resp. the orthogonal of $U$ w.r.t. the
pairing~\ref{eq:restrace}). Since $v_m=z_1^m$, the condition
${\rho(\sigma)U}=U$ is equivalent to:
$$
\res_{z=0}\tr\Big(\frac{z_1^m}{z_1}\psi_{\rho(\sigma)U}(z_1,t),
\frac{z_1^{1-m-p}}{z_1} \psi_U^*(z_1,s)\Big)dz\,=\,0$$
Bearing in
mind that $z_1^p=z$, together with the properties of the trace map, the
claim follows.

For the non-ramified case, recall that $U^\perp\in\gr^{-m}(V)$ for
$U\in\gr^m(V)$. Then, the condition ${\rho(\sigma)U}=U$ is written
as:
$$
\res_{z=0}\tr\Big(\frac{v_m}{z_\centerdot}\psi_{\rho(\sigma)U}(z_\centerdot,t),
\frac{v_{-m}}{z_\centerdot}
\psi_U^*(z_\centerdot,s)\Big)dz\,=\,0$$ The formula follows from
   the fact that $v_m v_{-m}=1$.
\end{pf}

Observe that the action considered in~\ref{eq:actiongroupongrv},
preserves the determinant bundle. To see this, it suffices to note
that for $\sigma\in G$  the commutative diagram:
   {\small
$$\xymatrix{ \L \ar[r] \ar[d]_{\rho(\sigma)} & V/V^+ \ar[d]^{\rho(\sigma)} \\
\rho(\sigma)\L \ar[r] & V/V^+}$$} implies the existence of an
isomorphism $\Det\simeq \rho(\sigma)^*\Det$.

This will allow us to find an explicit expression relating
the Baker-Akhiezer functions of $U$ and $\rho(\sigma)U$ in the
following general setting.

Assume that $V$ is the $\C((z))$-vector space
$V_1\times\dots\times V_r$ where $V_i=\C((z^{1/e_i}))$ and $e_i$
is a positive integer and that the action can be described as
follows. There exist natural numbers $\{k_1,\dots , k_l\}$ with
$m:=k_1+\dots+ k_l\leq r$ and a set $\{\xi_{m+1},\dots,\xi_r\}$
with $\xi_i$ a primitive $e_i$-th root of $1$ in $\C$ such that:
\begin{itemize}
\item $\rho(\sigma)$ acts as a cyclic permutation on the products:
   {\small $$\begin{aligned}
   & V_1\times\dots\times V_{k_1} \\
   &\quad\dots \\
   & V_{k_1+\dots+k_{l-1}+1}\times\dots\times V_{k_1+\dots+k_l}
   \end{aligned}$$}
\item
$\rho(\sigma)$ leaves  the factor $V_j$ (for $m<j\leq r$) stable
and maps $f(z^{1/e_j}))$ to $f(\xi_j z^{1/e_j})$.
\end{itemize}

\begin{prop}\label{prop:bafunctionsigmaU}
The  Baker-Akhiezer function of $\rho(\sigma)U$ is: {\small
$$\begin{aligned}
      \psi_{\rho(\sigma)U}({z_\centerdot},t)\,=\,
\psi_U\Big(\rho(\sigma){z_\centerdot}, &
    t^{(2)},\dots, t^{(k_1)},t^{(1)}, \dots \\
    & t^{(k_1+\dots+k_{l-1}+2)},\dots , t^{(k_1+\dots+k_{l})},
t^{(k_1+\dots+k_{l-1}+1)}, \\
& \xi_{m+1}^{-1}t_1^{(m+1)},\xi_{m+1}^{-2}t_2^{(m+1)}, \dots,
\xi_{r}^{-1}t_1^{(r)},\xi_{r}^{-1}t_2^{(r)},\dots\Big)
\end{aligned}
$$}
\end{prop}

\begin{pf}
It suffices to see that the automorphism $\rho(\sigma)$ gives
rise to a natural automorphism on the local Jacobian of the formal
spectral cover:
$$\rho(\sigma)^*\colon \j(\w C_V)\,\longrightarrow \, \j(\w C_V)$$
whose expression in terms of the underlying rings is: {\small $$
\begin{aligned}
& \rho(\sigma)^*(t_i^{(1)})=t_i^{(2)} \\ &\quad \dots \\ &
\rho(\sigma)^*(t_i^{(k_1)})=t_i^{(1)} \\  &\quad \dots \\ &
\rho(\sigma)^*(t_i^{(k_1+\dots+k_{l-1}+1)})=t_i^{(k_1+\dots+k_{l-1}+2)}
\\  &\quad \dots \\ &
\rho(\sigma)^*(t_i^{(k_1+\dots+k_{l})})=t_i^{(k_1+\dots+k_{l-1}+1)}
\\   &
      \rho(\sigma)^*(t_i^{(j)})=\xi_j^{-i}t_i^{(j)}\qquad\text{for }m<j\leq r
\end{aligned}
$$}
\end{pf}

\subsection{Grassmannian of isotropic subspaces}\hfill

For later purposes it is convenient to introduce a second
subscheme of the Grassmannian. Under our hypothesis one can
consider a skewsymmetric $p$-form on $V$. To this end, recall that
$V$ is a $p$-dimensional $\C((z))$-vector space; hence, mapping
$z_1\wedge\dots\wedge z_p$ (resp. $1\wedge z_1\wedge\dots \wedge
z_1^{p-1}$)
      to $dz$ we obtain an identification of $\wedge^p_{\C((z))} V$ with
$\C((z))dz$ for the non-ramified case (resp. ramified case).
Let us define a skewsymmetric form on $V$ by:
$$
\begin{aligned}
V\times\overset{p}\dots \times V \,& \longrightarrow\, \C \\
(f_1,\dots, f_p)\,& \longmapsto \, \res_{z=0}(f_1\wedge
\dots\wedge f_p)
\end{aligned}
$$
Let $\wedge^p_k W$ be the $p$-th exterior algebra of a $k$-vector
space $W$. For a $\C$-subspace $U$ of $V$, let us denote by $\wedge^p
U\subseteq \C((z))dz$ the image of the map:
      \beq
      \label{eq:lambdaUtoC((z))}
      \xymatrix{
      \wedge^p_{\C} U \ar@{^(->}[r] &  \wedge^p_{\C} V
       \ar@{->>}[r] & \wedge^p_{\C((z))} V
       \ar[r]^{\sim} & \C((z))dz}
\end{equation}

      A subspace $U$ of $V$ will be called isotropic if the
restriction of the skewlinear form to $U\times\dots\times U$ is
zero or, what amounts to the same, the map: $$\res_{z=0}\;\colon
\wedge^p U \longrightarrow \C $$ vanishes.

Consider the following subfunctor of the Grassmannian consisting
of the isotropic subspaces of $V$; that is:
      {\small \beq
      \label{eq:prymgrass}
S\,\rightsquigarrow\, \underline\pgrv(S)\,:=\,
\left\{\begin{gathered}
      U\in \underline\grv(S)\text{ such that the} \\
      \text{map }\wedge^p U \overset{\res}\longrightarrow \o_S\text{ is zero}
      \end{gathered}
      \right\}
      \end{equation}}

\begin{prop}
The functor $\underline\pgrv$ is representable by a closed
subscheme $\pgrv$ of $\grv$, which will be called the Grassmannian of
isotropic subspaces.
\end{prop}

\begin{pf}
Let $\o$ denote the sheaf of $\grv$ and let $\L\subset V\hat\otimes
\o$ denote the universal submodule. Then, the composition of the
multilinear form and the residue yields a morphism of locally free
$\o$-sheaves:
$$\L\otimes_{\o}\overset{p}\dots\otimes_{\o}\L\,\longrightarrow\,
\o$$ Since $\underline\pgrv$ consists of the points where this map
vanishes, the conclusion follows.
\end{pf}

Let us write down equations for this subscheme of $\grv$.

\begin{lem}\label{lem:generatorslambdaU}
Let $U$ be a point of $\gr^m(V)$ and let $\psi_U(z_\centerdot,t)$ be
its Baker-Akhiezer function. Let: {\small $$ \big(\psi_U^{1}(z,t),
\dots, \psi_U^{p}(z,t)\big)$$}
      be the coordinates of
$\frac{v_m}{z_\centerdot}\psi_U(z_\centerdot,t)$ w.r.t. the
distinguished basis of $V$ as
      a $\C((z))$-vector space.

Then, the $\C$-subspace $\wedge^p U$ of $\C((z))dz$ is generated
by the series of $\C((z))$ of the type: {\small
$$
\left| \matrix \psi_U^{1}(z,t_1) & \dots & \psi_U^{p}(z,t_1)
\\
\vdots &  & \vdots
\\
\psi_U^{1}(z,t_p) & \dots & \psi_U^{p}(z,t_p)
\endmatrix\right|
$$}
where $t_1,\dots, t_p$ denotes $p$ independent sets of variables
$t$.
\end{lem}

\begin{pf}
Since $\frac{v_m}{z_\centerdot}\psi_U(z_\centerdot,t)$ is a
generating function of $U$, it follows that $\wedge^p U$ is
generated by evaluating at certain values of $t_1,\dots , t_p$ the
image of
$\frac{v_m}{z_\centerdot}\psi_U(z_\centerdot,t_1)\wedge\dots\wedge
\frac{v_m}{z_\centerdot}\psi_U(z_\centerdot,t_p)$ by the
map~(\ref{eq:lambdaUtoC((z))}). To conclude, it suffices to express
the Baker-Akhiezer function of $U$ in terms of the distinguished
basis of $V$.
\end{pf}

\begin{thm}\label{thm:residueforpgrv}
Let $U$ be a closed point of $\gr^m(V)$. The point $U$ belongs to
$\pgrv$  if and only if the following identity holds: {\small \beq
\label{eq:bkpgeneral} \res_{z=0} \left| \matrix \psi_U^{1}(z,t_1)
& \dots & \psi_U^{p}(z,t_1)
\\
\vdots &  & \vdots
\\
\psi_U^{1}(z,t_p) & \dots & \psi_U^{p}(z,t_p)
\endmatrix\right| dz
\,=\, 0 \end{equation}}
\end{thm}

\begin{pf}
This is trivial from the previous lemma.
\end{pf}

\begin{rem}
It is worth noticing that the residue condition of the previous
theorem can be translated into a set of differential equations. To
this end, one should proceed as in Theorem~5.13 of~\cite{MP3}.
\end{rem}

\begin{rem}
If $p=2$ and $V=\C((z_1))$ (ramified case),  then the index $0$
connected component of $\pgrv$ coincides with the subscheme of
$\grv$ consisting of those subspaces that are maximal totally
isotropic:
{\small $$ \{U \in\gr^0(V) \text{ such that }
\res_{z_1=0}f(z_1)g(-z_1)\frac{dz_1}{z_1^2}=0\text{ for all
}f,g\in U\}$$}
and equation~(\ref{eq:bkpgeneral}) is
equivalent to the BKP hierarchy. That is, our approach generalizes
those of~\cite{PlPrym,Sh2}. Observe that the slight difference
between equation~(\ref{eq:bkpgeneral}) and the bilinear equation
of Lemma~4 in \cite{Sh2} stems from the fact that Shiota works
with the connected component of index $1$ while we do so with that of
index $0$.
\end{rem}

The relationship between the formal Prym variety and the
Grassmannian of isotropic subspaces is unveiled in the following:

\begin{thm}\label{th:accion pi}
The largest subgroup of $\Gamma_V$ acting on $\pgrv$ is $\Pi$.
\end{thm}

\begin{pf}
Observe that an element $g\in \Gamma_V$ acts on $\grv$, mapping a
point $U$ to $g\cdot U$, and that the homothety defined by $g$ is
an automorphism of $V$ as a $\C((z))$-vector space. It follows
that: \beq \label{eq:lambdanormU} \wedge^p (g\cdot U)\,=\, \nm(g)
\cdot \wedge^p U
\end{equation}

Then, it is trivial that $\Pi$ should act on the $\pgrv$ because the
norm of its elements belongs to $\C$.

Let us prove the converse.
Let us consider an element $g\in\Gamma_V$ acting on $\pgrv$. Since
the norm takes values in  $\C((z))$, it suffices to check that the
coefficient of $z^n$ in $\nm(g)$ is zero for all $n\neq 0$.

Given an integer $n\neq 0$, let  $U_n$ be a subspace of $V$ of the
type: {\small
$$ <\! z^{-n-1}e_1,\dots, e_p\!>\oplus z^N V^-\qquad N\in \Z
$$} where $\{e_1,\dots, e_p\}$ is the distinguished basis of $V$ and
$N\in \Z$ is such that $ U_n$
belongs to the $\pgrv$ (this condition is attained for small
enough $N$).

The construction of $U_n$ implies that
$z^{-n-1}e_1\wedge\dots\wedge e_p\in \wedge^p U_n$. Since $g\cdot
U_n\in \pgrv$, it follows by identity~(\ref{eq:lambdanormU})
that the coefficient of $z^n$ in $\nm(g)$ must be $0$ for all
$n\neq0$.
\end{pf}

\section{Prym varieties}

\subsection{Preliminaries}\hfill

Let $\pi\colon Y\to X$ be a cyclic covering of degree $p$ ($p$
being a prime number) between  smooth irreducible projective
curves over $\C$ of genus $g,\ \bar g$, respectively. Assume that
$\pi$ is given by an automorphism $\sigma_Y$ of $Y$ of order $p$;
in particular, $Y/<\!\sigma_Y\!>=X$. Let $R_\pi\subset Y$ denote
the ramification divisor of $\pi$  and consider the line bundle on
$X$ of degree $1/2\deg R_\pi$  given by $\delta=\wedge
\pi_*\o(R_\pi)$.

Let $J_d(Y)$ (resp. $J_d(X)$) denote the Jacobian of $Y$ (resp.
$X$) of degree $d$ and let $\nm\colon J_d(Y)\to J_d(X)$ be the
Albanese morphism induced by $\pi$. Note that $\wedge\pi_*L\simeq
\nm(L)\otimes \delta^*$ for every line bundle $L$ on $Y$. In
particular, for a line bundle $L\in J_d(Y)$ of degree:
$$d\,=\,(g-1)-(p-2)(\bar g-1)$$
   it follows that $\deg
(\wedge\pi_*L)=2(\bar g-1)$. Then, for this value of $d$, we
define the Prym variety of degree $d$ associated with the covering
$\pi$ as:
          $$ P_d(Y,\sigma_Y)\, :=\, \{ L\in J_d(Y)\, :\,
\wedge\pi_*L\simeq \omega_{X}\} \subset J_d(Y)$$
where
$\omega_{X}$ is the canonical sheaf on $X$. In other words:
$$L\in
P_d(Y,\sigma_Y)\quad\iff\quad \nm(L)\simeq \omega_{X}\otimes
\delta$$

We should be noted that the functor of points of $P_d(Y,\sigma_Y)$ is
given by
       $$P_d(Y,\sigma_Y)^\bullet (S)=\left\{ \begin{gathered} \, [L]\in
J_d(Y)^\bullet(S)\colon
\wedge (\pi\times\id)_*(L)\simeq q_1^*\omega_X\otimes q_2^*N \\
\text{ for some
line bundle $N$ on }S\end{gathered}\right\}$$
      for each $\C$-scheme $S$, where $[L]$ denotes the equivalence
class of a line bundle $L$ on $Y\times S\to S$ and $q_1,q_2$ are the
natural projections of $X\times S$.

\begin{rem}
It follows that the Prym variety of degree $d$ is the translation
of $\ker\nm \subset J_0(Y)$ by an arbitrary $L_0\in P_d(Y)$.
Observe that the classical Prym variety associated with $\pi$ is the
connected component of $\o_Y$ in $\ker\nm$.
\end{rem}

In this section the following data will be fixed $(Y,\sigma_Y,\bar
y,t_{\bar y})$ where:
\begin{enumerate}
\item $Y$ is a smooth irreducible projective curve over $\C$ of genus
$g$;
\item $\sigma_Y$ is an order $p$ automorphism  of $Y$;
\item $\bar y$ is an orbit  of $\sigma_Y$; and,
\item $t_{\bar y}\colon \w\o_{Y,\bar y}\overset\sim\to V^+$ is a
formal parameter along $\bar y$ that is equivariant with respect
to the actions of $\sigma_Y$ and $\sigma$.
\end{enumerate}

Let $\pic^{\infty}(Y,\bar y)$ be the scheme whose functor of
points is given by:
$$S \,\rightsquigarrow\, \pic^{\infty} (Y,\bar y)^\bullet(S)\, =\,
\left\{  \begin{gathered} \text{pairs $(L,\phi)$ where }
       L  \text{ is a line bundle} \\
\text{on  $Y\times S$ and }\phi\colon \w L_{\bar y\times
S}\overset\sim\to \o_S\w \otimes V^+\end{gathered}\right\} $$
We
denote by $\pic^{\infty}_d(Y,\bar y)$ the subscheme of
$\pic^\infty (Y,\bar y)$ consisting of those pairs $(L,\phi)$,
where the line bundle $L$ has degree $d$.

The Krichever map is the injective morphism of functors given by:
$$\begin{aligned} K\colon \pic^\infty (Y,\bar y)^\bullet (S)
&\,\longrightarrow \,\grv^\bullet (S)
\\ (L,\phi) & \,\longmapsto \, \phi\big(\limil {n\geq 0} (p_Y)_*L(n \cdot
(\bar y\times S))\big)\end{aligned}$$
where $p_Y\colon Y\times S\to S$ is
the natural projection (for a more detailed study of $\pic^\infty (Y,\bar
y)$ and $K$, see ~\cite{A}).

   Consider the
point in $X$ given by  $x=\pi(\bar y)$. Since $t_{\bar y}$ is
equivariant, it gives rise to an isomorphism between the
subrings of invariant elements; that is, a formal parameter
$t_x\colon \w\o_{X,x}\overset\sim\to \C[[z]]$  at $x$ such that
the diagram:
$$\xymatrix { \w\o_{Y,\bar y}\ar[r]_{\sim}^{t_{\bar y}} &  V^+ \\
\w\o_{X,x} \ar[r]_{\sim}^{t_x} \ar@{_{(}->}[u] & \C[[z]]
\ar@{_{(}->}[u]}$$ is commutative. Moreover, note that $t_x$
induces a map:
$$dt_x\colon (\omega_X)^{\w{\!}}_x \,\overset\sim\to\,
     \Omega_{\widehat \o_{X,x}/\C} \,\overset\sim\to\,
     \Omega_{\C[[z]]/\C} \,\overset\sim\to\,
     \C[[z]]dz$$

If we are given a pair $(L,\phi)$ in $\pic^{\infty}_d (Y,\bar y)$,
then the determinant of $\phi$ yields an isomorphism: {\small
$$(\det\phi)\colon ( \wedge  \pi_*L)^{\w{\!}}_x\, \overset\sim\to \,
\wedge_{\C[[z]]} V^+\, \overset\sim\to \, \C[[z]] z_1\wedge \dots
\wedge z_p$$}

     When $\omega_X\simeq \wedge  \pi_*L$, we have an isomorphism (up to a
non-zero constant): $(\omega_X)^{\w{\!}}_x\overset\sim\to (\wedge
\pi_*L)^{\w{\!}}_x$. We require these isomorphism to be
compatible:

\begin{defn}
If $\omega_X\simeq \wedge  \pi_*L$, the formal parameter $t_{\bar
y}$ and the formal trivialization $\phi$ of $L$ along $\bar y$ are
said to be compatible if the diagram: {\small
$$\xymatrix{(\omega_X)^{\w{\!}}_x \ar[r]^\sim
\ar[d]^\simeq_{dt_x} & (\wedge  \pi_*L)^{\w{\!}}_x \ar[d]_\simeq^{\det\phi}
\\
\C[[z]]dz \ar[r]^-{\sim} & \C[[z]] z_1\wedge \dots \wedge z_p}
$$}
commutes (up to a non-zero constant), where the bottom arrow maps
$dz$ to $z_1\wedge \dots \wedge z_p$ and $t_x$ is the formal
parameter at $x$ defined by $t_{\bar y}$.
\end{defn}

Similarly, one defines compatibility for $S$-valued points.
Indeed, let $t_{\bar y}$ be a formal trivialization of $\o_Y$
along $\bar y$ and let $(L,\phi)$ be a point in $\pic^\infty(Y,\bar
y)^\bullet(S)$ such that $L\in P_d(Y,\sigma_Y)^\bullet(S)$. Then,
$t_{\bar y}$ and $\phi$ are said to be compatible, if the diagram:
{\small
$$\xymatrix{(\omega_X)^{\w{\!}}_x\w\otimes_\C \o_S \ar[r]^\sim
\ar[d]^\simeq_{dt_x\otimes 1} & (\wedge
(\pi\times\id)_*L)^{\w{\!}}_{x\times S} \ar[d]_\simeq^{\det\phi}
\\
\o_S[[z]]dz \ar[r]^-{\sim} & \o_S[[z]] z_1\wedge \dots \wedge z_p}
$$}
commutes (up to an element of $H^0(S,\o_S^*)$).

\begin{defn}
The functor $\fprym^{\infty} (Y,\sigma_Y,\bar y)$ associated with
the data $(Y,\sigma_Y,\bar y, t_{\bar y})$ is the subfunctor of
$\pic^{\infty}_d (Y,\bar y)^\bullet$ defined by: {\small $$S
\,\rightsquigarrow\, \fprym^{\infty} (Y,\sigma_Y,\bar y)(S)\, :=\,
\left\{ \begin{gathered} (L,\phi)\, \in\, \pic^{\infty}_d (Y,\bar
y)^\bullet(S)\,
       \text{ such} \\
       \text{that }
\wedge (\pi\times \id)_*L\simeq q_1^*\omega_X
\\
\text{and $t_{\bar y}$ and $\phi$ are compatible}
\end{gathered}\right\} $$}
\end{defn}

\begin{thm}
The functor $\fprym^{\infty} (Y,\sigma_Y,\bar y)$ is representable
by a closed subscheme $\prym^{\infty} (Y,\sigma_Y,\bar y)$ of
$\pic^{\infty}_d (Y,\bar y)$.
\end{thm}

\begin{pf}
It suffices to check that the injective morphism of functors:  {\small $$
\begin{aligned}
\fprym^{\infty} (Y,\sigma_Y,\bar y) &\hookrightarrow
\pic^{\infty}_d (Y,\bar y)^\bullet  \underset {J_d(Y)^\bullet}\times
P_d(Y,\sigma_Y)^\bullet  \\ (L,\phi) &\mapsto \big( (L,\phi)\, , \,
[L]\big)\end{aligned}
$$}
(where $[L]\in J_d(Y)^\bullet(S)$ denotes the equivalence class of
$L$) is a closed immersion; in other words, that given an
$S$-valued point:
$$\big( (L,\phi)\, , \, [L] \big)\in
\pic^{\infty}_d (Y,\bar y)^\bullet (S) \underset {J_d(Y)^\bullet
(S)}\times P_d(Y,\sigma_Y)^\bullet (S)$$
the condition  that this point belongs to $\fprym^{\infty}
(Y,\sigma_Y,\bar y)(S)$ is a closed condition on $S$.

One checks that $\wedge (\pi\times\id)_*L\simeq q_1^*\omega_X\otimes
q_2^*N$ for some line bundle $N$ on $S$. Therefore, considering
formal completions along ${x}\times S$, the formal parameter $t_x$
at $x$ induces an isomorphism (up to a non-zero constant):
$$\begin{aligned} ( \wedge  (\pi\times \id)_*L)^{\w{\;}}_{{x}\times S}\,
       &\simeq
\, ( (\omega_{X})^{\w{\;}}_x\, \widehat \otimes_\C\, \o_S)\,
\widehat\otimes_{\o_S}\, N\, \simeq \\
& \simeq \, (\Omega_{\C[[z]]/\C}\, \widehat \otimes_\C\, \o_S)\,
\widehat\otimes_{\o_S}\, N\, \simeq \, \o_S[[z]]dz\,
\widehat\otimes_{\o_S}\, N\end{aligned}$$

On the other hand, the determinant of $\phi$ yields an
isomorphism: {\small $$(\det\phi)\colon ( \wedge  (\pi\times
\id)_*L)^{\w{\!}}_{{x}\times S}\, \overset\sim\to \,
\wedge_{\o_S[[z]]} (V^+\, \widehat\otimes_\C\, \o_S)\,
\overset\sim\to \, \o_S[[z]] z_1\wedge \dots \wedge z_p$$}

Comparing both expressions, it follows that  $N\simeq \o_S$.  Fixing an
isomorphism $q_1^*\omega_X\simeq\wedge(\pi\times \id)_*L $, the condition
that $t_{\bar y}$ and $\phi$ are compatible is equivalent to saying that
the
$\o_S[[z]]$-module isomorphism:
{\small
$$\o_S[[z]]dz\overset{(dt_x\otimes 1)^{-1}}{\underset{\sim}\longrightarrow}
(\omega_{X})^{\w{\;}}_x\, \underset\C{\widehat \otimes}\,
\o_S\underset{\sim}\longrightarrow ( \wedge  (\pi\times
\id)_*L)^{\w{\;}}_{{x}\times S}\overset{\det
\phi}{\underset{\sim}\longrightarrow} \o_S[[z]]z_1\wedge \dots\wedge z_p$$}
sends $dz\mapsto \lambda\cdot z_1\wedge \dots \wedge z_p$ with
$\lambda\in H^0(S,\o_S)$, which is a closed condition. Hence, the
statement follows.
\end{pf}

\begin{prop}\label{prop:picgrvi=prym}
One has the following cartesian diagram:
$$\xymatrix {\pic^\infty_d (Y,\bar y)\ar@{^{(}->}[r]^-{K} & \grv  \\
\prym^\infty (Y,\sigma_Y,\bar y)\ar@{^{(}->}[r]^-{K} \ar@{_{(}->}[u]
& \pgrv \ar@{_{(}->}[u]}$$
\end{prop}

\begin{pf}
If suffices to check the claim for geometric points.
Let $(L,\phi)\in \prym^\infty (Y,\sigma_Y,\bar y)$ and
$U=K(L,\phi)=\phi\big(H^0(Y-\bar y, L)\big)\subset
V$.

Since $\wedge\pi_*L\simeq \omega_X$ and $t_{\bar y}$ and $\phi$
are compatible, we have the commutativity (up to a constant) of:
{\small
$$\xymatrix{(\omega_X)^{\w{\!}}_x \ar[r]^\sim
\ar[d]^\simeq_{dt_x} & (\wedge  \pi_*L)^{\w{\!}}_x \ar[d]_\simeq^{\det\phi}
\\
\C[[z]]dz \ar[r]^-{\sim} & \C[[z]] z_1\wedge \dots \wedge z_p}
$$}

This diagram implies that the composition map: {\small
$$\wedge_\C^pH^0(X-x,\pi_*L) \twoheadrightarrow
H^0(X-x,\wedge \pi_*L)\simeq H^0(X-x,\omega_X)
\overset{dt_x}\hookrightarrow \C((z))dz$$}
coincides (up a
multiplicative constant) with
morphism~(\ref{eq:lambdaUtoC((z))}) given in the definition of
$\pgrv$. In particular, it turns out that: $$\wedge^p U\, =\,
dt_x\big(H^0(X-x,\omega_X)\big)$$ as subspaces of $\C((z))dz$.

Since $H^0(X-x,\omega_X)$ are meromorphic differentials with one
only pole, the  condition of the vanishing of the residue (at
$z=0$) is satisfied; that is, $U=K(L,\phi)=\phi(H^0(Y-\bar y, L))\in
\pgrv$.

     Conversely, let $U\in\pgrv$ such that there exists $(L,\phi)\in
\pic_d^{\infty}(Y, \bar y)$ with $K(L,\phi)=U$. The above
arguments show that: $$\wedge^p U\,
=\,(\det\phi)\big(H^0(X-x,\wedge\pi_*L)\big)$$ as subspaces of
$\C((z))dz$.

The non-degenerated pairing:
{\small $$\begin{aligned} \C((z))\times \C((z))dz &\longrightarrow \C \\
(f,\omega )&\longmapsto \res_{z=0}(f\cdot\omega)\end{aligned}$$}
induces isomorphisms:
{\small $$\begin{aligned} R\colon \gr\C((z))&\longrightarrow
\gr\big(
\C((z))dz\big)\\
R'\colon\gr\big(\C((z))dz\big)&\longrightarrow \gr
\C((z))\end{aligned}$$}
such that $R\circ R'=\id_{\gr\C((z))dz}$ and $R'\circ R=\id_{\gr \C((z))}$.

Let $A_{\wedge^pU}=t_x(H^0(X-x,\o_X))\subset \C((z))$. The
condition $U\in\pgrv$ implies that the map:
$$\wedge^p U \,=\,(\det \phi)\big(
H^0(X-x,\wedge\pi_*L)\big)\,\subset\,
\C((z))dz\,\overset{\res}\longrightarrow \, \C$$
vanishes. Then, it follows that
$A_{\wedge^pU}\subseteq R'(\wedge^p U)$ because $f\cdot
\wedge^pU\subset \wedge^pU$ for all $f\in A_{\wedge^pU}$.

Using the Krichever construction for the triple $(X,x,t_x)$, one
checks that the subspace $R'(\wedge^p U)\subset \C((z))$ is
attached to the geometric datum $(\omega_X\otimes
(\wedge\pi_*L)^{-1}, dt_x\otimes ((\det\phi)^*)^{-1})$, where
$(\det\phi)^*$ denotes the transposed map of $\det\phi$. From this, one
deduces that $A_{\wedge^pU}$ and $R'(\wedge^p U)$ lie in the same
connected component of $\gr\C((z))$. Therefore,
$A_{\wedge^pU}= R'(\wedge^p U)$.

Bearing in mind that $R\circ R'=\id_{\gr\C((z))dz}$, it holds that:
$$dt_x(H^0(X-x,\omega_X))=R(A_{\wedge^pU})=\wedge^p
U=(\det\phi)(H^0(X-x,\wedge\pi_*L))$$
By the injectivity of the Krichever morphism, this identity implies that
$\omega_X\simeq \wedge\pi_*L$ and that $t_{\bar y}$ and $\phi$ are
compatible. That is, $(L,\phi)\in \prym^\infty (Y,\sigma_Y,\bar y)$.
\end{pf}

\begin{thm}
The formal group $\Pi$ acts on the image of the Krichever map $
K\colon \prym^\infty (Y,\sigma_Y,\bar y)\hookrightarrow \pgrv$.
\end{thm}

\begin{pf}
Bearing in mind that $\Pi\subset \Gamma_V$ acts on the image of
the Krichever map $ K\colon \pic^\infty_d (Y,\bar
y)\hookrightarrow \grv$, the proof follows from Theorem~\ref{th:accion pi}
      and
the above proposition.
\end{pf}

\subsection{Geometric characterization of
Pryms}\label{subsec:geometriccharacterization}\hfill

In this part we shall show an analogue for Prym varieties of
Mulase's characterization of Jacobians varieties (\cite{Mul}).
Roughly speaking, he shows that Jacobian varieties are precisely
the finite dimensional orbits of a certain group acting on a
quotient Grassmannian. Our approach will work within the frameset
of formal schemes (\cite{EGA},~I.\S10) and follow~\cite{PlPrym}
closely. This characterization is strongly related to
that of Shiota,
which claims that Jacobians are the finite
dimensional solutions of the KP hierarchy (\cite{Sh,Sh2}).

The action of $\Gamma_V$ on $V$ by homotheties gives rise to an
action:    {\small
$$
\begin{aligned}
\mu\colon \Pi\times \grv^p &\,\longrightarrow\, \grv^p \\
(g,U_{\centerdot})&\,\longmapsto \, g\cdot U_{\centerdot}
\end{aligned}
$$}
where $\grv^p:=\grv\times\overset{p}\dots\times\grv$,
$U_{\centerdot}:=(U_1,\dots, U_p)\in\grv^p$ and $g\cdot
U_{\centerdot}:=(g\cdot U_1,\dots,g\cdot U_p)$.

On the other hand,  the following map: {\small $$\begin{aligned}
\grv\,& \hookrightarrow\,
\grv^p \\
U\, &\mapsto U_{\sigma}:=(U,\sigma(U),\dots,\sigma^{p-1}(U))
\end{aligned}$$}
turns out to be a closed immersion.

\begin{defn}
Given a closed point $U\in\grv$ we define the morphism $\mu_U$ by:
{\small $$\begin{aligned} \mu_U\colon\Pi \,& \longrightarrow\,
\grv^p \\
g\, &\mapsto g\cdot U_{\sigma}=(g \cdot U,g \cdot
\sigma(U),\dots,g \cdot \sigma^{p-1}(U))
\end{aligned}$$}
\end{defn}

Our goal consists in characterizing those points of $\grv$ coming
from geometric data in terms of the above morphism. To this end,
we begin with some general properties.

\begin{lem}
Let $U$ be a closed point of $\grv$. The orbit of $\U_\sigma\in
\grv^p$ under $\mu$ is the  schematic image of $\mu_U$ which is a
formal scheme and will be denoted by $\Pi({U_\sigma})$.
\end{lem}

\begin{pf}
Following Remark~\ref{rem:jisdirectlim}, we write the group $\Pi$
as a direct limit of affine group schemes $\Pi^n$.

Let us denote by $\o$ the sheaf of rings of $\grv^p$, by
$\mu_{U,n}^* :\o\to \o_{\Pi^n}$ the morphism induced by
$\mu_U\vert_{\Pi^n}$, and by $I_n$ the ideal $\ker \mu_{U,n}^*$.
Finally, define $\Pi^n({U_\sigma})$ to be the schematic image of
the morphism: $$\mu_U\vert_{\Pi^n}\colon\Pi^n\,\longrightarrow\,
\grv^p$$ that is:
$$\Pi^n({U_\sigma})\, :=\, \sp\big(\o/I_n\big)$$

     For $m\geq n\geq 0$, consider the
morphisms $\phi^m_n$ and $\bar \phi^m_n$ defined by the following
commutative diagram: {\small $$ \xymatrix @R=8pt { & \o/I_m
\ar@{^(->}[r] \ar[dd]^{\phi^m_n} & \o_{\Pi^m} \ar@{->>}[dd]^{\bar
\phi^m_n}
\\
      \o \ar@{->>}[ru] \ar@{->>}[rd]
\\
& \o/I_n \ar@{^(->}[r]  & \o_{\Pi^n} }$$}

Since $\Pi({U_\sigma})=\cup_{n\geq 0} \Pi^n({U_\sigma})$, in order
to show that $\limil{n} \Pi^n({U_\sigma})$ exists in the category
of formal schemes it suffices  to show that  $\phi^m_n$ is
      surjective and that its kernel is a nilpotent ideal
(\cite{EGA},~I.\S10.6.3).

Bearing in mind that $I_k=\ker(\mu_{U,k}^*)$ and that
$\mu_{U,n}^*=\bar \phi_n^m\circ \mu_{U,m}^*$, the surjectivity of
$\phi^m_n$ follows easily. Finally, note that the ideal
$\ker(\phi^m_n)$ is contained in $\ker(\bar \phi^m_n)$, which is
nilpotent.
\end{pf}

\begin{lem}
The stabilizer of ${U_\sigma}$ is a closed subgroup of $\Pi$,
which will be denoted by $H_{U_\sigma}$, and there is a canonical
isomorphism of formal schemes:
$$\Pi/H_{U_\sigma}\,\simeq \, \Pi({U_\sigma})$$
(We understand by a closed subgroup of $\Pi$ a subgroup defined by a
closed ideal of $\o_\Pi$).
\end{lem}

\begin{pf}
Recall that $\Pi=\limi \Pi^n$, where $\Pi^n$ are affine group
schemes. For $n\geq 0$, consider the subfunctor $\underline
H^n_{U_\sigma}$ whose set of $R$-valued points ($R$ being a
$\C$-algebra) is:
$$\underline H^n_{U_\sigma}(R)\,:=\, \{ g\in\Pi^n(R) \,\vert\, g\cdot
U_{\sigma}=U_{\sigma}\}$$ Since $\underline H^n_{U_\sigma}$ is
defined by a closed condition, it is representable by a closed
subscheme $H^n_{U_\sigma}$ of $\Pi^n$, which is precisely the
stabilizer of ${U_\sigma}$ under the action of $\Pi^n$ on
$\grv^p$.

Note that the following cartesian diagram (for $m\geq n\geq 0$):
\begin{equation} \label{eq:cartesian} \xymatrix{ H^n_{U_\sigma} \ar@{^(->}[r]
\ar@{^(->}[d] & \Pi^n
\ar@{^(->}[d]\\
H^m_{U_\sigma} \ar@{^(->}[r] & \Pi^m} \end{equation}
     implies that the morphism $\o_{H^m_{U_\sigma}}\to \o_{H^n_{U_\sigma}}$ is
     surjective and its kernel is a nilpotent ideal. It follows that
$H_{U_\sigma}:=\limi H^n_{U_\sigma}$ is a formal scheme and a
closed subgroup of $\Pi$ that coincides with the stabilizer of
${U_\sigma}$.

Let us now prove  the second part. Bearing in mind that
$H_{U_\sigma}^n$ is the stabilizer soubgroup of the action of
$\Pi^n$ and that $\Pi^n({U_\sigma})$ is the schematic image of
$\mu_U\vert_{\Pi^n}$, one has that:
$$\Pi({U_\sigma})\,=\, \limi \Pi^n({U_\sigma})\,=\, \limi
\Pi^n/H^n_{U_\sigma}$$
We conclude if we can show that $ \limi \Pi^n/H^n_{U_\sigma} =
\Pi/H_{U_\sigma}$. However, this follows easily from the following
commutative diagram:
$$\xymatrix{
0 \ar[r] & H^n_{U_\sigma} \ar[r] \ar@{^(->}[d] & \Pi^n \ar[r]
\ar@{^(->}[d] & \Pi^n/H^n_{U_\sigma} \ar[r] \ar@{^(->}[d] & 0 \\
0 \ar[r] & H^m_{U_\sigma} \ar[r]  & \Pi^m \ar[r]&
\Pi^m/H^m_{U_\sigma} \ar[r] & 0 }
$$
where the injectivity of the third vertical arrow follows from
the fact that the diagram~\ref{eq:cartesian} is cartesian.
\end{pf}

For the sake of notation, it is convenient to consider the
following subgroup of $\Pi$:
$$\bar \Pi^+=\cases \G_m\times \Pi^+ \qquad
\hphantom{GGGGGGGt}\text{in the ramified case,}\\
(\G_m\times\overset{p}\dots\times\G_m)\times\Pi^+\qquad \text{in
the non-ramified case.}\endcases$$

\begin{lem}
The quotient $\Pi({U_\sigma})/\bar \Pi^+$ is a formal scheme
canonically isomorphic to $\Pi/H_{U_\sigma}+ \bar \Pi^+$. In
particular, its tangent space at $U_\sigma$ is isomorphic to
$T_1\Pi/(\ker d\mu_U+ T_1\bar \Pi^+)$ where:
$$d\mu_U\colon T_1\Pi\,\longrightarrow\, T_{U_\sigma}\grv^p$$
is the map induced by $\mu$ on the tangent spaces.
\end{lem}

\begin{pf}
In order to show that $\Pi({U_\sigma})/\bar \Pi^+$ is a formal
scheme, it suffices to check that it is a direct limit of affine
schemes in the conditions of~\cite{EGA},~I.\S10.6.3.  Since
$\Pi^n({U_\sigma})=\Pi^n/H^n_{U_\sigma}$, one has that
$\Pi^n({U_\sigma})/\bar \Pi^+=\Pi^n/H^n_{U_\sigma}+\bar \Pi^+$.

Finally, the commutative diagram:
$$\xymatrix{
0 \ar[r] & H^n_{U_\sigma}+\bar \Pi^+ \ar[r] \ar@{^(->}[d] & \Pi^n
\ar[r]
\ar@{^(->}[d] & \Pi^n/H^n_{U_\sigma}+\bar \Pi^+ \ar[r] \ar@{^(->}[d] & 0 \\
0 \ar[r] & H^m_{U_\sigma}+\bar \Pi^+ \ar[r]  & \Pi^m \ar[r]&
\Pi^m/H^m_{U_\sigma}+\bar \Pi^+ \ar[r] & 0 }
$$
implies that the direct limit of $\Pi^n/H^n_{U_\sigma}+\bar \Pi^+$
exists as a formal scheme; it will be denoted by
$\Pi/H_{U_\sigma}+\bar \Pi^+$.

The claim about the tangent space is now straightforward since
$T_1 H_{U_\sigma}\simeq \ker d\mu_U$.
\end{pf}

\begin{rem}
Observe that $d\mu_U$ is explicitly given by: {\small
\beq\label{eq:kerdmuU}
\begin{aligned}
d\mu_U\colon T_1\Pi \,& \longrightarrow\,  T_{U_\sigma}\grv^p\simeq
\prod_{i=0}^{p-1}\hom(\sigma^i(U),V/\sigma^i(U)) \\
g\, &\longmapsto \, \,\big\{\xymatrix@C=14pt{\sigma^i(U)
\ar@{^(->}[r] & V \ar[r]^{\cdot g} & V \ar@{->>}[r] &
V/\sigma^i(U)}\big\}_{i=0,\dots,p-1}
\end{aligned}
\end{equation}}
and that $T_1\Pi=\{g\in V\vert \tr(g)\in \C\}$ because the norm
$\nm:\Gamma_V \to \Gamma$ induces the trace at their tangent
spaces. Therefore, one has that: {\small $$\begin{aligned} \ker
d\mu_U \,& =\, \{ g\in V \,\vert \, \tr(g)\in\C \text{ and }
g\cdot \sigma^i(U)\subseteq \sigma^i(U) \text{ for all } i\} \,=\\
& =\, \{ g\in V \,\vert \, \tr(g)\in\C \text{ and }
\sigma^i(g)\cdot U\subseteq U \text{ for all } i\}
\end{aligned}
$$}
\end{rem}

\begin{thm}\label{thm:algebraizableifffinitedimension}
Let $U$ be a closed point of $\grv$. Then, the following
conditions are equivalent:
\begin{enumerate}
\item
$\Pi({U_\sigma})/\bar \Pi^+$ is algebraizable,
\item
$\dim_\C T_{U_\sigma}(\Pi({U_\sigma})/\bar \Pi^+) <\infty$,
\end{enumerate}
Further, if one of these conditions holds true, then the ring of
$\Pi({U_\sigma})/\bar \Pi^+$ is isomorphic to a ring of series in
finitely many variables.
\end{thm}

\begin{pf}
$1\implies 2$. Let us write $\Pi({U_\sigma})/\bar \Pi^+$ as $\sf
A$ and let us note that $A$ is an admissible linearly topologized ring.
Recall that algebraizable (\cite{Hart}~II.9.3.2) means that the
formal scheme is isomorphic to the completion of a noetherian
scheme along a closed subscheme. Since the completion of a
noetherian ring with respect to an ideal is noetherian, one
concludes that $A$ is noetherian.

Let $J:=\limp J_n$ be a definition ideal of $A$. Since $A$ is
noetherian and $A/J\simeq k$, it follows from \cite{EGA}~0.7.2.6
that $\left( J/J^2\right)^*$ is a finite dimensional vector space.

We conclude from the following inclusion: {\small
$$\aligned T_{U_\sigma}\big(\Pi({U_\sigma})/\bar \Pi^+\big)\,&=\,
\hom_{\text{for-sch}}\big(\sp(k[\epsilon]/(\epsilon^2)),\sf A\big)\, =\\
&=\, \hom \Sb \text{topological} \\ \text{$k$-algebras} \endSb
(A,k[\epsilon]/(\epsilon^2)) \,\subseteq \\
&\subseteq \,
\hom_{\text{$k$-algebras}}(A,k[\epsilon]/(\epsilon^2))
\,\subseteq\, \left( J/J^2\right)^*
\endaligned$$}

$2\implies 1$. To prove the claim, it will be enough to show that
$A$ is isomorphic to a ring of power series in finitely many
variables, $\C[[u_1,\dots, u_d]]$, endowed with the $(u_1,\dots,
u_n)$-adic topology.

Recall that for all $n>\!>0$ there is a surjection
     of vector spaces:
$$\xymatrix{ T_1 \Pi^n/\bar \Pi^+\ar@{->>}[r] &
T_{U_\sigma}\big(\Pi^n({U_\sigma})/\bar \Pi^+\big)}$$
Let us introduce the
following notation: $\sp(A_n):=\Pi^n({U_\sigma})/\bar \Pi^+$, $\sf
A:=\limi\sp(A_n)$ and  $\sp(\o_n):=\Pi^n/\bar \Pi^+$.  Bearing in mind
that $A_n$ and $\o_n$ are artinian rings, one has that:
$${\frak m}_{A_n}/{\frak m}^2_{A_n}\,\hookrightarrow\, {\frak m}_{\o_n}/{\frak
m}^2_{\o_n}\qquad\forall n>\!>0$$ where ${\frak m}_{A_n}$ (resp.
${\frak m}_{\o_n}$) denotes the maximal ideal of $A_n$ (resp.
$\o_n$).

   Recall from Remark~\ref{rem:jisdirectlim} that $\o_n$ is
of the type $\C[u_{n'},\dots, u_1]/I_{n'}$, with $n'$ depending on
$n$, and where $(u_{n'},\dots, u_1)^n\subseteq I_{n'}$.

On the other hand, the hypothesis claims that the vector space:
$$ T_{U_\sigma}\big(\Pi({U_\sigma})/\bar \Pi^+\big)\,=\, \limi
T_{U_\sigma}\big(\Pi^n({U_\sigma})/\bar \Pi^+\big) $$ is of finite
dimension, say $d$. That is:
$$\dim \big(
T_{U_\sigma}\big(\Pi^n({U_\sigma})/\bar \Pi^+\big)\big) \,= \,
d\qquad \forall n>\!> 0$$

   Since $T_{U_\sigma} \sp(A_n)\simeq
({\frak m}_{A_n}/{\frak m}_{A_n}^2)^*$ and the maps:
$$\begin{gathered}
J:=\limp {\frak m}_{A_n}\,\overset{\pi_n}\longrightarrow\,  {\frak m}_{A_n} \\
{\frak m}_{A_m}\,\longrightarrow\,  {\frak
m}_{A_n}\end{gathered}$$ are surjective, there exist elements
$\bar v_1,\dots,\bar v_d\in J$ such that:
$$<\!\{\pi_n(\bar v_1),\dots, \pi_n(\bar v_d)\}\!>={\frak
m}_{A_n}/{\frak m}_{A_n}^2\qquad\text{ for all }n>\!>0$$ By
Nakayama's lemma one has epimorphisms:
$$\aligned p_n:k[v_1,\dots,v_d]&\to A_n\\
v_i\quad&\mapsto\pi_n(\bar v_i)\endaligned\qquad \forall n>\!>0$$
compatible with the  surjections $A_m\to A_n$ for $m\geq n>\!>0$.

Summing up, we have obtained a commutative diagram:
$$\xymatrix{
<\! v_1,\dots, v_d\!> \ar@{=}[d] \ar[r]^\sim & {\frak
m}_{A_n}/{\frak m}^2_{A_n} \ar@{^(->}[r] \ar@{=}[d] & <\!
u_1,\dots,
u_{n'}\!> \ar@{->>}[d] \\
<\! v_1,\dots, v_d\!>  \ar[r]^\sim & {\frak m}_{A_m}/{\frak
m}^2_{A_m} \ar@{^(->}[r]  & <\! u_1,\dots, u_{m'}\!> }$$ for
$n\geq m>\!>0$.

Up to a change of coordinates, it can be assumed that the image of
$v_i$ is $u_i$ for all $1\leq i\leq d$. Looking at the
corresponding morphisms of rings, we find the following
commutative diagram:
$$\xymatrix{
\C[v_1,\dots, v_d]/ (v_1,\dots, v_d)^{n} \ar@{->>}[d] \ar@{->>}[r]
& {A_n} \ar@{^(->}[r] \ar@{->>}[d] & \C[u_1,\dots,
u_{n'}]/I_{n'}  \ar@{->>}[d] \\
\C[v_1,\dots, v_d]/ (v_1,\dots, v_d)^{m}  \ar@{->>}[r] & {A_m}
\ar@{^(->}[r]  & \C[u_1,\dots, u_{m'}]/I_{m'} }$$
Taking inverse
limits, we have that the following morphism of linearly
topologized rings:
$$\xymatrix{
\C[[v_1,\dots, v_d]] \ar[r] & A=\limp A_n \ar@{^(->}[r]
&\C\{\{u_1,u_2\dots\}\} }$$ maps $v_i$ to $u_i$. In particular, this
implies that:
$$A\,= \, \C[[u_1,\dots, u_d]] $$
with the $(u_1,\dots, u_d)$-adic topology and the conclusion
follows.
\end{pf}

\begin{thm}\label{thm:algebraizableiffgeometric}
Let $U$ be a closed point of $\grv$. Then, the following
conditions are equivalent:
\begin{enumerate}
\item
$\Pi({U_\sigma})/\bar \Pi^+$ is algebraizable,
\item
there exist data $(Y,\sigma_Y, \bar y,t_{\bar y})$ and a line
bundle with a formal trivialization $(L,\phi)$  such that $\bar y$
is an orbit of $\sigma_Y$ and:
$$U=(t_{\bar y}\circ\phi)(H^0(Y-\bar y,L))$$
\end{enumerate}
\end{thm}

\begin{pf}
$1\implies 2$. For $U\in\grv$, consider the vector spaces: {\small
$$\begin{aligned} \ker d\mu_U\,&
=\, \ker \big( T_1\Pi\to T_{U_\sigma}\grv^p\big)\,=\\
& =\,  \{ g\in V \,\vert \, \tr(g)\in\C \text{ and } g\cdot
\sigma^i(U)\subseteq \sigma^i(U) \text{ for all } i\} \\
   B\,:&  =\,
\ker
\big(T_1\Gamma_V \to T_{U_\sigma}\grv^p\big)\,=\\
& =\,
     \{ g\in V \,\vert \,
g\cdot \sigma^i(U)\subseteq \sigma^i(U) \text{ for all } i
\}\end{aligned}$$}

Note that the exact sequence of formal group schemes: {\small
$$
\xymatrix{ 0 \ar[r] & \Pi  \ar[r] & \Gamma_V \ar[r]^{\nm}  &
\Gamma \ar[r] & 0}$$}
     induces the exact sequence of vector spaces:
{\small
$$
\xymatrix{ 0 \ar[r] & T_1\Pi  \ar[r] & T_1\Gamma_V= V \ar[r]^{\tr}
& T_1\Gamma = \C((z)) \ar[r] & 0}$$} This sequence allows us to
consider the following diagram of $\C$-vector spaces: {\small
$$
\xymatrix{ 0 \ar[r] & \ker d\mu_U \ar[r] \ar[d] & B \ar[r]^{\tr}
\ar[d] & \tr(B) \ar[r] \ar[d] & 0
\\
0 \ar[r] & T_1\Pi/T_1\bar \Pi^+ \ar[r] & V/V^+ \ar[r] & \C((z))/\C[[z]]
\ar[r] & 0 }$$}

  Snake's Lemma yields the following exact sequence:
{\small \begin{multline}\notag
0  \to  \ker d\mu_U\cap  T_1\bar \Pi^+ \to  B\cap V^+ \to \tr(B)\cap
\C[[z]] \to \\
    \hphantom{mm}\to T_1\Pi/(\ker d\mu_U+ T_1\bar \Pi^+) \to V/(B+V^+) \to \\
\to \C((z))/(\tr(B)+\C[[z]]) \to 0
\end{multline}}

    From the expressions of $\ker d\mu_U$ and $B$, one easily checks
that the first two terms of the above sequence are finite
dimensional vector spaces.

Observe that the subalgebra of $V$ generated by $ \ker d\mu_U$ is
a subalgebra of $B$ of finite codimension (because $p$ is a prime
number). From hypothesis (1) and
Theorem~\ref{thm:algebraizableifffinitedimension}, we know that
$T_1\Pi/(\ker d\mu_U+ T_1\bar \Pi^+)$ is finite dimensional, so
$V/(B+V^+)$ is finite dimensional too. Summing up, we have
obtained a pair $(U,B)\in\grv\times \grv$ such that $B$ is a
sub-$\C$-algebra of $V$ and $U$ is a $B$-module.

Standard arguments of the Krichever construction show that the
pair $(U,B)$ correspond to a curve $Y$, a divisor $\bar y$, a
formal parameter along the divisor $t_{\bar y}:\widehat\o_{Y,\bar
y}\simeq V^+$, a line bundle $L$ and a formal trivialization
$\phi:\widehat L_{\bar y}\simeq \widehat\o_{Y,\bar y}$
(\cite{MP3} Theorem~5.3, \cite{Mul2, Ines}). Further, it is also known
that, under the Krichever correspondence, we have:
$$\begin{gathered}
U\,=\, (t_{\bar y}\circ\phi)(H^0(Y-\bar y, L)) \\
B\,=\, t_{\bar y}(H^0(Y-\bar y, \o_Y)) \end{gathered}$$

The above discussion also shows that $B$ carries an action of
$\sigma$ and, therefore, the curve $Y$ is canonically endowed with
an order $p$ automorphism $\sigma_Y$.

$2\implies 1$. Given  $(Y,\sigma_Y,\bar y,t_{\bar y})$ and
$(L,\phi)$  as in the statement, one considers the quotient $\pi:
Y\to X:=Y/<\!\sigma_Y\!>$. Let $B$ be the $\C$-algebra $t_{\bar
y}(H^0(Y-\bar y,\o_Y)) \,\in \, \grv$ which is endowed with an
action of $\sigma$ induced by $\sigma_Y$. Since $\sigma$ acts on
$B$ and $\tr(b)=\sum_{i=0}^{p-1}\sigma^i(b)$, it follows that:
$$\tr(B)\,\subseteq\, B$$
and (from Theorem~5.10 of~\cite{MP3}):
$$\tr(B) \,=\, t_x(H^0(X-\pi(\bar y),\o_X)) \,\in \, \gr\C((z))$$

Summing up, the maps $B\to V/V^+$ and $\tr(B)\to \C((z))/\C[[z]]$
have finite dimensional kernels and cokernels. Then, using the
long exact sequence of the first part of the proof, we conclude
that $T_1\Pi/(\ker d\mu_U+ T_1\bar \Pi^+)$ is finite dimensional or,
what amounts to the same, $\Pi({U_\sigma})/\bar \Pi^+$ is algebraizable
(Theorem~\ref{thm:algebraizableifffinitedimension}).
\end{pf}

\begin{cor}\label{cor:tangenttoorbitiskerneloftrace}
If the conditions of Theorem~\ref{thm:algebraizableiffgeometric}
are satisfied, then there is an isomorphism:
$$T_{U_\sigma}(\Pi({U_\sigma})/\bar \Pi^+)\,\simeq\, T_{\o_Y}P_d(Y,\sigma_Y)$$
\end{cor}

\begin{pf}
Recall that in the proof of
Theorem~\ref{thm:algebraizableiffgeometric}, we obtained that
$\tr(B)\subseteq B$, which implies that the connecting morphism of
the long exact sequence of that proof is zero; that is, the
sequence splits into two exact short sequences. The second one
reads as follows:
$$0\to T_{U_\sigma}(\Pi({U_\sigma})/\bar \Pi^+)\to
     H^1(Y,\o_Y)\overset{\tr}\to H^1(X,\o_X)
     \to 0$$
which implies  the  claim.
\end{pf}

\begin{thm}\label{thm:characterizationofpryms}
Let $U$ be a closed point of $\grv$ such that  $\Pi({U_\sigma})/\bar
\Pi^+$ is algebraizable, or equivalently, there exist data
$(Y,\sigma_Y,\bar y,t_{\bar y})$ and  $(L,\phi)$ with $U=(t_{\bar
y}\circ\phi)(H^0(Y-\bar y,L))$. Then, it holds that:
$$U\,\in\,\pgrv \ \iff\ L\,\in\,P_d(Y,\sigma_Y) \text{ and } \phi
\text{ and } t_{\bar y} \text{ are compatible.}$$
\end{thm}

\begin{pf}
This follows from Proposition~\ref{prop:picgrvi=prym} and
Theorem~\ref{thm:algebraizableiffgeometric}.
\end{pf}

\begin{rem}
The corresponding Prym flows have been computed by several authors
in terms of Lax pairs (\cite{AB,K}) or in terms of the Heisenberg
algebras (\cite{DM}).
\end{rem}

\section{Equations of the moduli space of curves with automorphisms}

The aim of this section is to explicitly describe the moduli space
of pointed curves with a non-trivial automorphism group as a
subscheme of the Sato Grassmannian. Let us recall that given a
smooth irreducible projective curve $Y$ of genus $g$ over $\C$ and
$p$ a prime number such that $p$ divides the order of $\aut (Y)$,
one has that $p\leq 2g+1$. Then, the problem of characterizing the
moduli space of curves of genus $g$ with a non-trivial automorphism
group is reduced to the characterization of the moduli space of
curves with an order $p$ automorphism, where  $p$ is a prime number
and $p\leq 2g+1$. Henceforth, $p$ will denote a prime number.

Let $Y$ be a smooth irreducible projective curve over $\C$, let
$\sigma_Y\in \aut (Y)$ be an order $p$ automorphism and let
$X=Y/<\!\sigma_Y\!>$ be the quotient under the action of
$<\!\sigma_Y\!>$. The projection $\pi\colon Y\to X$ is a cyclic
covering of degree $p$. Given  a geometric point $x\in X$, there
are two posibilities for its fibre:
\begin{enumerate}
\item[(a)] {\sl Ramified case:} $x$ is a ramification point of $\pi$;
that is, the fibre  contains
only one point $y\in Y$:{\small  $$\pi^{-1}(x)\,=\, p\cdot y$$}
\item[(b)] {\sl Non-ramified case:} $x$ is not a ramification point of
$\pi$; that is,  the fibre contains $p$ pairwise different points
$y_1,\ldots , y_p\in Y$: {\small
      $$\pi^{-1}(x)=y_1+\dots +y_p$$}
\end{enumerate}

We shall study the moduli of pairs $(Y,\sigma_Y)$
distinguishing these two cases: ramified and non-ramified. It is
therefore convenient to use different notations for the two cases
considered in section \S2:
\begin{enumerate}
\item[(a)] {\sl Ramified case:} $V_{\mathrm R}= \C((z^{1/p}))=\C((z_1))$ and
$V_{\mathrm R}^+=\C[[z_1]]$.
\item[(b)] {\sl Non-ramified case:} $V_{{\mathrm NR}}=\C((z_1))\times
\dots\times \C((z_p))$ and $V_{{\mathrm NR}}^+=\C[[z_1]]\times
\dots\times \C[[z_p]]$.
\end{enumerate}
Following the notations of \S 3.A, we shall denote by
$\gr(V_{\mathrm R})^\sigma$ and $\gr(V_{{\mathrm NR}})^\sigma$ the
closed subschemes whose points are the invariant points under the
action of $<\!\sigma\!>$.

\begin{defn}
Let $\M^{\infty,1}$ be the closed subscheme of $\gr(V_{\mathrm
R})$ representing the moduli of integral complete curves with one
marked smooth point and a formal parameter at this point (see
\cite{MP}, Theorem 6.5).

Let $\M^{\infty,p}$ be the closed subscheme of $\gr(V_{{\mathrm
NR}})$ representing the mo\-du\-li of reduced complete curves with
$p$ marked pairwise distinct smooth points and formal parameters
at these points (see \cite{MP3}, Theorem 5.3).
\end{defn}

\subsection{The ramified case and the Krichever map}\hfill

\begin{defn}
The moduli functor $\fm^{\infty}(p,{\mathrm R})$ of curves with an
automorphism of order $p$ and a formal parameter at a fixed point is the
contravariant functor on the category of
$\C$-schemes:
$$
\fm^{\infty}(p,{\mathrm R})\colon
\left\{{\scriptsize\begin{gathered}\text{category of}\\
\C-\text{schemes}\end{gathered}}\right\}\,
\rightsquigarrow\,\left\{{\scriptsize\begin{gathered}\text{category}\\\text{of
sets}\end{gathered}}\right\}
$$
that associates the set of equivalence
classes  of data $\{Y,\sigma_Y, y,t_y\}$ with a $\C$-scheme $S$, where:
\begin{enumerate}
\item
$p_Y:Y\to S$ is a proper and flat morphism whose fibres are
geometrically integral curves.
\item $\sigma_Y$ is a non-trivial automorphism of the curve $Y\to S$
of order $p$ such
that the quotient $X:=Y/<\!\sigma_Y\!>$ exists and  $p_X:X\to
S$ is a proper and flat morphism whose fibres are geometrically integral
curves.
\item $y:S\to Y$ is a smooth section of $p_Y$ which is a Cartier
divisor of $Y$ invariant under the action of $\sigma_Y$ and such
that, for each closed point $s\in S$, $y(s)$ is a smooth point of
the fibre $Y_s=p_Y^{-1}(s)$.
\item There exists a formal parameter $t_y$ along the section $y$:
      $$t_y\colon \w\o_{Y,y(S)}\overset\sim\longrightarrow
V_{\mathrm R}^+\,\w\otimes\, \o_S$$
that is compatible with
respect to the actions of $\sigma_Y$ and $\sigma$.
\item
$\{Y,\sigma_Y,y,t_y\}$ and $\{Y',\sigma'_Y,y',t_{y'}\}$ are said to be
equivalent when there is a commutative diagram of  $S$-schemes: {\small
$$\xymatrix{ Y \ar[d]_{\pi} \ar[r]^{\sim} & Y' \ar[d]^{\pi'} \\
X \ar[r]^{\sim} & X'}$$} compatible with all the data, where $\pi$, $\pi'$
are the natural projections.
\end{enumerate}
\end{defn}

Let us define the Krichever map for $\fm^\infty(p,{\mathrm R})$ as
the morphism of functors:
      $$K\colon \fm^\infty(p,{\mathrm R}) \longrightarrow \gr(V_{\mathrm
R})^ \bullet$$ such that to each $\{ Y,\sigma_Y,y,t_y\}\in
\fm^\infty(p,{\mathrm R})(S)$ one attaches the submodule:
      $$K(Y,\sigma_Y,y,t_y)\,=\,t_y\,\big( \,\limil {n\geq
0}(p_Y)_*\o_Y(n\cdot y(S))\big) \,\subset\, V_{\mathrm R}\,
\w\otimes_\C\, \o_S$$

\begin{thm}\label{th:ramificado}
The functor $\fm^\infty(p,{\mathrm R})$ is represented by a closed
subscheme $\M^\infty(p,{\mathrm R})$ of $\gr(V_{\mathrm R})$ that
coincides with $\M^{\infty,1}\cap \gr(V_{\mathrm R})^\sigma$.
\end{thm}

\begin{pf}
One has only to prove the inclusion:
      $$(\M^{\infty,1})^\bullet(S) \cap
(\gr(V_{\mathrm R})^\sigma)^\bullet(S)\,\subseteq \, K\,\big(\,
\fm^\infty(p,{\mathrm R})(S)\,\big)$$ for any $\C$-scheme $S$. If
$U$ is a point of $(\M^{\infty,1})^\bullet(S) \cap (\gr(V_{\mathrm
R})^\sigma)^\bullet(S)$ representing a geometric datum
$(Y,y,t_y)$, then the automorphism:
$$\sigma\otimes\id\colon
V_{\mathrm R}\,\w\otimes_\C\, \o_S\to V_{\mathrm
R}\,\w\otimes_\C\, \o_S$$ leaves the $\o_S$-module
$U\subset V_{\mathrm R}\,\w\otimes_\C\, \o_S$ stable.

Since these functors are sheaves, we can assume that $S=\sp (B)$
is an affine scheme. Therefore, $U\subset B((z_1))=V_{\mathrm
R}\,\w\otimes_\C\, B$ is a sub-$B$-algebra endowed with the
natural filtration induced by the degree of $z_1^{-1}$. The curve
$Y$ can be reconstructed from $U$ as $Y=\proj_B\U$, where $\U$ is
the graded algebra associated with the above filtration
(see~\cite{MP}, Theorem 6.4, \cite{Ines}). Since the automorphism
$\sigma\otimes\id$ induces an automorphism of the $B$-algebra $U$
compatible with the filtration, then it induces an automorphism
$\wtilde\sigma\colon \U\to \U$ of graded $B$-algebras.
     From this fact, one easily concludes that we can construct from $U$ a
geometric datum $(Y,\sigma_Y,y,t_Y)$ defining a point of
$K\,\big(\, \fm^\infty(p,{\mathrm R})(S)\,\big)$.
\end{pf}

In order to be able to write down equations for $\M^\infty
(p,{\mathrm R})$, note that the Euler-Poincar\'{e} characteristic
of the curve must be fixed. Thus, let us denote by
$\M^\infty_\chi(p,{\mathrm R})$ the subscheme of $\M^\infty
(p,{\mathrm R})$ representing curves whose structure sheaf has
Euler-Poincar\'{e} characteristic $\chi$. We have that
$\M^\infty_\chi(p,{\mathrm R})\subset \gr^{\chi}(V_{\mathrm R})$.

Then, we can prove the following theorem:

\begin{thm}
The closed subscheme $\M^\infty_\chi(p,{\mathrm R})$ of
$\gr^{\chi}(V_{\mathrm R})$ is defined by the following set of
equations:
\begin{enumerate}
\item $
\res_{z=0}\tr\Big(\frac1{z_1}\psi_U(\xi
z_1,\xi^{-1}t_1,\xi^{-2}t_2,\ldots )
\psi_U^*(z_1,t)\Big)\frac{dz}z=0$.
\item
$\res_{z=0}\tr\Big(z_1^{\chi-2}\psi_U(z_1,t)\psi_U(z_1,t')\psi^*_U(z_1,t'')\Big)
\frac{dz}z\,=\,0$.
\item
$\res_{z=0}\tr\Big(\frac1{z_1^{\chi}}\psi^*_U(z_1,t)\Big)
\frac{dz}z\,=\,0$.
\end{enumerate}
\end{thm}

\begin{pf}
Bearing in mind the expression for the Baker-Akhiezer function
given in Proposition~\ref{prop:bafunctionsigmaU} and the
characterization given in
Theorem~\ref{thm:characterizationgrsigma}, one has that the
equations defining the subscheme $\gr^{\chi}(V_{\mathrm
R})^\sigma$ are given by (1).

Recall that  $U\in \gr^{\chi}(V)$ lies in $\M^{\infty,1}$ if and
only if $U\cdot U\subseteq U$ and $1\in U$ (\cite{MP},
Theorem~6.4). Note that $z_1^{\chi-1}\psi_U(z_1,t)$ (resp.
$\frac1{z_1^{p+\chi}}\psi^*_U(z_1,t)$) is a generating function of
$U$ (resp. $U^\perp$). In terms of the pairing defined by
formula~\ref{eq:restrace}, these two conditions read:
$$\res_{z=0}\tr\Big(z_1^{2(\chi-1)}\psi_U(z_1,t)\psi_U(z_1,t'),\frac1{z_1^{p+\chi}}\psi^*_U(z_1,t'')\Big)
dz\,=\,0$$ and
$$\res_{z=0}\tr\Big(1, \frac1{z_1^{p+\chi}}\psi^*_U(z_1,t)\Big)
dz\,=\,0$$ respectively. From the properties of the trace and
since $z_1^p=z$, these equations are equivalent to (2) and (3).
\end{pf}

\subsection{The non-ramified case and the Krichever map}\hfill

Study of the non-ramified case is quite similar to the
previous one.

\begin{defn}
The moduli functor $\fm^{\infty}(p,{\mathrm NR})$ of curves with
an automorphism of order $p$ and  formal parameters at the points
of a non-ramified orbit is the contravariant functor on the
category of $\C$-schemes:
$$
\fm^{\infty}(p,{\mathrm NR})\colon
\left\{{\scriptsize\begin{gathered}\text{category of}\\
\C-\text{schemes}\end{gathered}}\right\}\,
\rightsquigarrow\,\left\{{\scriptsize\begin{gathered}
\text{category}\\\text{of
sets}\end{gathered}}\right\}
$$
that associates the set of equivalence
classes  of data $\{Y,\sigma_Y,\bar y,t_{\bar y}\}$ with a $\C$-scheme $S$,
where:
\begin{enumerate}
\item
$p_Y:Y\to S$ is a proper and flat morphism whose fibres are
geometrically reduced curves.
\item $\sigma_Y$ is a non-trivial automorphism of the curve $Y\to S$
of order $p$ such
that the quotient $X:=Y/<\!\sigma_Y\!>$ exists and
$p_X:X\to S$ is a proper and flat morphism whose fibres are geometrically
integral curves.
\item There exist $p$ disjoint smooth sections $y_j\colon S\to Y$
$(j=1,\ldots ,p)$ of $p_Y$ which are Cartier divisors of $Y$ and such that
$\bar y=\{ y_1,\ldots ,y_p\}$ is an orbit under the action of
$<\!\sigma_Y\!>$.
\item For each closed point $s\in S$, $y_j(s)$ is a smooth point of the
fibre $p_Y^{-1}(s)=Y_s$ such that each irreducible component of
$Y_s$ contains at least a point $y_j(s)$ for a certain $j=1,\ldots
,p$.
\item There exists formal parameters $t_{y_j}$ along the sections $y_j$
such that they induce an  isomorphism:
      $$t_{\bar y}\colon \w\o_{Y,\bar y(S)}\overset\sim\longrightarrow
V_{{\mathrm NR}}^+\,\w\otimes\, \o_S$$ which is compatible with respect to
the actions of $\sigma_Y$ and $\sigma$.
\item
$\{Y,\sigma_Y,\bar y,t_{\bar y}\}$ and $\{Y',\sigma'_Y,\bar
y',t_{\bar y'}\}$ are said to be equivalent when there is a commutative
diagram of
$S$-schemes: {\small
$$\xymatrix{ Y \ar[d]_{\pi} \ar[r]^{\sim} & Y' \ar[d]^{\pi'} \\
X \ar[r]^{\sim} & X'}$$} compatible with all the data where $\pi$, $\pi'$
are the natural projections.
\end{enumerate}
\end{defn}

We define the Krichever map for $\fm^\infty(p,{\mathrm NR})$ as
the morphism of functors:
$$K\colon \fm^\infty(p,{\mathrm NR}) \longrightarrow \gr(V_{{\mathrm
NR}})^ \bullet$$
such that for each $\{ Y,\sigma_Y,\bar y,t_{\bar y}\}\in
\fm^\infty(p,{\mathrm NR})(S)$ one has:
      $$K(Y,\sigma_Y,\bar y,t_{\bar y})\,=\,t_{\bar y}\,\big( \,\limil {n\geq
0}(p_Y)_*\o_Y(n\cdot \bar y(S))\big)\, \subset \, V_{{\mathrm
NR}}\, \w\otimes_\C\, \o_S$$

\begin{thm}\label{th:noramificado}
The functor $\fm^\infty(p,{\mathrm NR})$ is represented by the
closed subscheme
      $\M^{\infty,p}\cap
\gr(V_{{\mathrm NR}})^\sigma\,=\,\M^\infty(p,{\mathrm NR})$
      of $\gr(V_{{\mathrm NR}})$.
\end{thm}

\begin{pf}
This is similar to the proof of Theorem~\ref{th:ramificado}.
\end{pf}

It is worth pointing out that in our definition of
$\M^\infty(p,{\mathrm NR})$, non-connected curves have been
considered. Since $p$ is a prime number, we have two posibilities:
either the curve $Y$ is connected or has $p$ connected components
such that each of them is isomorphic to the quotient curve $X$.

Our aim is to write down equations for $\M^\infty(p,{\mathrm NR})$.
As in the previous case, we must restrict ourselves to a
given Euler-Poincar\'{e} characteristic. Thus, let
$\M_{\chi}^\infty(p,{\mathrm NR})$ denote the subscheme of
$\M^\infty(p,{\mathrm NR})$ representing curves such that the
Euler-Poincar\'{e} characteristic of the structure sheaf is $\chi$.

If the Euler-Poincar\'{e} characteristic is positive, then there are
only   two cases; namely, $\chi=1$ and $\chi=p$. The first case
corresponds to the projective line endowed with an automorphism of
order $p$. The second one consists of the disjoint union of $p$
copies of the projective line where the automorphism is given by a
cyclic permutation of the copies. Therefore, the non-trivial cases
of curves with automorphisms appear when $\chi\leq 0$.

\begin{thm}\label{thm:equationsNR}
Let $\chi\leq 0$. The subscheme $\M_\chi^\infty(p,{\mathrm NR})$
of $\gr^{\chi}(V_{{\mathrm NR}})$ is defined by the following set
of equations:
\begin{enumerate}
\item
$\sum_{i=1}^p\res_{z=0}\psi_{\sigma(U)}^{(i)}(z,t)\cdot
\psi_U^{*(i)}(z,t')\,\frac{dz}{z^2}\, =\, 0$, or equivalently:
{\footnotesize
$$\res_{z=0}\psi_U^{(p)}(z,t)\cdot
\psi_U^{*(1)}(z,t')\frac{dz}{z^2} + \cdots +
\res_{z_p=0}\psi_U^{(1)}(z,t)\cdot
\psi_U^{*(p)}(z,t')\frac{dz}{z^2} = 0$$}
\item
   {\small $ \begin{aligned} & \sum_{i=1}^r \res_{z=0}\Big(\psi^{(i)}_U(z,t)
   \psi_U^{(i)}(z,t)\psi_U^{*(i)}(z,t)
   \Big)\frac{dz}{z^{q+4}}+
   \\ & \quad +
   \sum_{i=r+1}^p \res_{z=0}\Big(\psi^{(i)}_U(z,t)
   \psi_U^{(i)}(z,t)\psi_U^{*(i)}(z,t)
   \Big)\frac{dz}{z^{q+3}}\,=\,0 \end{aligned} $}
\item
    {\small $ \sum_{i=1}^r \res_{z=0}\psi_U^{*(i)}(z,t)
   z^q  dz+
   \sum_{i=r+1}^p \res_{z=0}\psi_U^{*(i)}(z,t)
   z^{q-1} dz \,=\,0 $}
\end{enumerate}
where $q, r$ are given by $-\chi=q\cdot p+r$ with $0\leq r< p$.
\end{thm}

\begin{pf}
Recalling Proposition~\ref{prop:bafunctionsigmaU} and
Theorem~\ref{thm:characterizationgrsigma}, one sees that the
equations defining the subscheme $\gr^{\chi}(V_{\mathrm
R})^\sigma$ are given by (1).

It remains to compute the equations of $\M_\chi^{\infty,p}$. These
equations are derived from the known fact that a point
$U\in\gr^\chi(V)$ lies in $\M_\chi^{\infty,p}$ if and only if it
is a $\C$-algebra with unity.

Since $\frac{v_{\chi}}{z_\centerdot}\psi_U(z_\centerdot,t)$ (resp.
$\frac{v_{-\chi}}{z_\centerdot}\psi_U(z_\centerdot,t)$) is a
generating function of $U$ (resp. $U^\perp$), the point $U$ is an
$\C$-algebra with unity precisely when: {\small
$$ \left\{\begin{aligned} &
   \res_{z=0}\tr\Big( \frac{v_{\chi}}{z_\centerdot}\psi_U(z_\centerdot,t)
   \cdot \frac{v_{\chi}}{z_\centerdot}\psi_U(z_\centerdot,t)
   \cdot \frac{v_{-\chi}}{z_\centerdot}\psi_U^*(z_\centerdot,t)
   \Big)dz\,=\, 0
   \\
   &
   \res_{z=0}\tr\Big( \frac{v_{-\chi}}{z_\centerdot}\psi_U^*(z_\centerdot,t)
   \Big)dz\,=\, 0
   \end{aligned}\right.$$}
If $-\chi=qp+r$ with $0\leq r< p$, then $v_{-\chi}=z_1^{q+1}\dots
z_r^{q+1}z_{r+1}^q\dots z_p^q$ and $v_{\chi}:=v_{-\chi}^{-1}$.
Recalling that $z_\centerdot=(z_1,\dots,z_p)$ and that: {\small
$$\tr\big((f_1(z_1),\dots,f_p(z_p)\big)\,=\,\sum_{i=1}^p
f_i(z)$$}one concludes.
\end{pf}

    The  points associated with non-connected curves define a closed
    subscheme of
$\M_{\chi}^\infty(p,{\mathrm NR})$ (by~\cite{MP3},
Proposition~3.20). Let $\M_{\chi}^\infty(p,{\mathrm NR})^0$ denote
the open subscheme parametrizing points such that the
corresponding curve is connected.

Bearing in mind that the Euler-Poincar\'{e} characteristic of a
non-connected curve with an order $p$ automorphism such that the quotient
curve is connected is a multiple of
$p$, it follows that $\M_{\chi}^\infty(p,{\mathrm NR})^0$
coincides with $\M_{\chi}^\infty(p,{\mathrm NR})$ provided that
$\chi\neq \dot p$.

Let us characterize $\M_{\chi}^\infty(p,{\mathrm NR})^0$ for
$\chi=\dot p$ and $\chi\leq 0$. Consider the elements: {\small
$$e_i\,:=\, (0,\ldots,\overset{(i)}1,\ldots ,0)\,\in \, V_{{\mathrm
NR}}\,=\, \C((z_1))\times \cdots \times \C((z_p))$$}
Then, a point
$U\in \M_{\chi}^\infty(p,{\mathrm NR})$ is associated with a
non-connected curve if and only if $e_i\in U$ for some $i\in \{
1,\ldots ,p\}$. This condition is equivalent to saying that one term
of the third sum of the statement of Theorem~\ref{thm:equationsNR}
does vanish; that is: {\small
$$\res_{z=0}\psi_U^{*(i)}(z,t)
   z^q  dz\,=\, 0\qquad \text{ for some } i\in\{1,\dots,p\}$$}
   (observe that $r=0$ since $\chi=\dot p$).

Since $\sigma$ acts transitively on the subset $\{ e_1,
\ldots,e_p\}$ and $U$ is invariant under $\sigma$, it follows that
if $e_i\in U$ for some $i$, then $e_i\in U$ for all $i\in \{
1,\ldots ,p\}$. Therefore, we deduce the
following corollary:

\begin{cor} The open subscheme $\M_{\chi}^\infty(p,{\mathrm NR})^0$  of
$\M_{\chi}^\infty(p,{\mathrm NR})$ is defined by the condition: {\small
$$\res_{z=0}\psi_U^{*(i)}(z,t)
   z^q  dz\,\neq\, 0\qquad \text{ for some } i\in\{1,\dots,p\}$$}
for $-\chi=q\cdot p\geq 0$, and  $\M_{\chi}^\infty(p,{\mathrm
NR})^0=\M_{\chi}^\infty(p,{\mathrm NR})$ for $\chi\neq \dot p$.
\end{cor}



\end{document}